%% file: Discontinuity_yildiz_etds.tex
\title{Discontinuity of Topological Entropy for the Lozi Maps}
\author{
         Izzet Burak Yildiz \\
                }
\date{}

\documentclass[12pt]{article}
\usepackage{amsmath , amsthm , commath, epsf}
\usepackage{amssymb}
\usepackage{amsfonts}
\usepackage{pause}
\usepackage{graphicx}
\usepackage{url}

\include{macros}

\begin{document}
\maketitle

\begin{abstract}
Recently, Buzzi \cite{Buzzi} showed in the compact case that the entropy map $f\rightarrow$ $h_{top}(f)$ is lower semi-continuous for all piecewise affine surface homeomorphisms. We prove that topological entropy for the Lozi maps can jump from zero to a value above $0.1203$ as one crosses a particular parameter and hence it is not upper semi-continuous in general. Moreover, our results can be extended to a small neighborhood of this parameter showing the jump in the entropy occurs along a line segment in the parameter space.
\end{abstract}

\section{Introduction}
There have been some recent developments in the study of piecewise affine surface homeomorphisms. In the compact case, Buzzi proved that under the assumption of positive topological entropy, there are finitely many ergodic measures maximizing the entropy \cite{Buzzi}. He also showed that topological entropy is lower semi-continuous for these maps. The following question was asked by Buzzi: \\

\textbf{Question:} \emph{Prove or disprove the upper semi-continuity of entropy for piecewise affine homeomorphisms of the plane.} \\

Our goal is to answer Buzzi's above question in the non-compact case by showing that topological entropy of the Lozi map is not upper semi-continuous at a given parameter. Moreover, our results can be extended to show that there is a line segment in the parameter space along which the topological entropy is not upper semi-continuous. \\

Let us start with a review of the subject:

\emph{Piecewise affine homeomorphisms:} Let $f : \mathbb{R}^n \to \mathbb{R}^n$ be a homeomorphism where $n \in \mathbb{Z}^{+}$. An \emph{affine subdivision} of $f$ is a finite collection $\mathcal{U}=\{U_1,\dots,U_N\}$ of pairwise disjoint non-empty open subsets of $\mathbb{R}^n$ such that their union is dense in $\mathbb{R}^n$ and $f|_{U_i} = A_i|_{U_i}$ for each $i=1,\dots,N$ where $A_i : \mathbb{R}^n \to \mathbb{R}^n$ is an invertible affine map. A \emph{piecewise affine homeomorphism} is a homeomorphism $f : \mathbb{R}^n \to \mathbb{R}^n$ for which there exists an affine subdivision.

\textbf{Example:} Lozi maps are piecewise affine homeomorphisms of the plane given by:
\[ \mathcal L = \mathcal L_{a,b} : \left( \begin{array}{ccc} x \\ y \end{array} \right) \mapsto \left( \begin{array}{ccc} 1-a|x|+by \\ x \end{array} \right) \textit{,} \hspace{4mm} a \textit{,} b \in \mathbb{R} \textit{,} \hspace{2mm} b \neq 0 .
\]
Note that $\mathcal{U}=\{U_1,U_2\}$ where $U_1=\{(x,y)\in \mathbb{R} \hspace{2mm}|\hspace{2mm} x>0\}$ and $U_2=\{(x,y)\in \mathbb{R} \hspace{2mm}|\hspace{2mm} x<0\}$. \\

Let us first review some of the results about continuity properties of entropy in different dimensions.
Throughout this paper, we will denote the topological entropy of a map $f$ by $h(f)$. \\

In one dimension, one can work with piecewise monotone functions. Let $I$ denote a compact interval of $\mathbb{R}$. A map $T: I \to I$ is called a \emph{piecewise monotone} function if there exists a partition of $I$ into finitely many subintervals on each of which the restriction of $T$ is continuous and strictly monotone. Two piecewise monotone maps $T_1$ and $T_2$ are said to be $\varepsilon$-close, if they have the same number of intervals of monotonicity and the graph of $T_2$ is contained in an $\varepsilon$-neighborhood of the graph of $T_1$ considered as subsets of $\mathbb{R}^2$. It was proved by Misiurewicz and Szlenk \cite{Mis1} that the entropy map $f \to h(f)$ is lower semi-continuous for piecewise monotone continuous maps. They also gave upper bounds for the jumps up of the entropy. For unimodal maps, entropy is continuous for all maps for which it is positive \cite{Mis5}. \\

There are also some continuity results in higher dimensions. Let $C^r(M^n)$ denote the set of $C^r$ self maps of an $n$-dimensional compact manifold. It is a classical result of Katok \cite{Katok2} that the entropy map is lower semi-continuous for $C^{1+\alpha}$ diffeomorphisms on compact surfaces. Yomdin \cite{Yomdin} and Newhouse \cite{Newhouse2} proved that entropy is upper semi-continuous in $C^{\infty}(M^n)$ for $n\geq1$. Combining these two results, one can get the continuity of entropy in $C^{\infty}(M^2)$. This result does not hold for homeomorphisms on surfaces \cite{Rees}. Also, Misiurewicz \cite{Mis4} constructed examples showing that entropy is not continuous in $C^{\infty}(M^n)$ for $n\geq4$ as well as examples \cite{Mis2} showing that entropy is not upper semi-continuous in $C^r(M^n)$ where $r<\infty$ and $n\geq2$. \\

For piecewise affine surface homeomorphisms, the following Katok-like theorem (see \cite{Katok1}) was given by Buzzi \cite{Buzzi}:

\begin{Theorem} Let $f : M \to M$ be a piecewise affine homeomorphism of a
compact affine surface. Let $S$ be the singularity locus of $M$, that is, the set of
points $x$ which have no neighborhood on which the restriction of $f
$ is affine.
For any $\varepsilon > 0$, there is a compact invariant set $K \subset M \setminus S$ such that
$h(f|K) > h(f) -\varepsilon$.
Moreover $f : K \to K$ is topologically conjugate to a subshift of finite type.
\end{Theorem}

The lower semi-continuity of the entropy in the compact case follows from the above theorem. This result may also hold in the non-compact case but it requires more work. The goal of this paper is to disprove the upper semi-continuity in the non-compact case by showing a jump up of the entropy in Lozi maps. Our results can be summarized as follows:

\begin{Theorem}[\bf{Main Theorem}] \label{A} In general, the topological entropy of the Lozi map does not depend continuously on the parameters:
There exists some $\epsilon_*>0$ such that for all $0<\epsilon_1<\epsilon_*$ and $|\epsilon_2|<\epsilon_*$,
\begin{itemize}
\item[(i)]  The topological entropies of the Lozi maps with $(a,b)=(1.4+\epsilon_2,0.4+\epsilon_2)$, $h(\mathcal L_{1.4+\epsilon_2,0.4+\epsilon_2})$, are zero.
\item[(ii)]  The topological entropies of the Lozi maps, $h(\mathcal L_{(1.4+\epsilon_1+\epsilon_2,0.4+\epsilon_2)})$, have a lower bound of $0.1203$.
\end{itemize}

\end{Theorem}
In other words, we show that the entropy is zero on the line segment $l=\{(1.4 + \epsilon_2, 0.4 +\epsilon_2): |\epsilon_2|<\epsilon_*\}$ and it is above $0.1203$ for the parameters immediately to the right of that segment.

\section{Topological Entropy}

Topological entropy is a quantitative measurement of how complicated a map is.

\begin{Definition} Let $f:X\to X$ be a continuous map on a compact metric space $(X,d)$ with a metric $d$. Two distinct points $x,y \in X$, $x\neq y$, are called \emph{$(n,\epsilon)$-separated} for a positive integer $n$ and $\epsilon >0$ if there is $m \in \{0,1,\ldots,n-1\}$, such that $d(f^m(x), f^m(y)) > \epsilon$. A set $U\subset X$ is called an \emph{ $(n,\epsilon)$-separated set} if every pair of distinct points $x,y \in U$, $x\neq y$, is $(n,\epsilon)$-separated. \\

Let $r(n,\epsilon,f)$ be the maximum cardinality of an $(n,\epsilon)$-separated set $U\subset X$. By compactness, this number is always finite. Define $h(\epsilon, f)$ = $\displaystyle\limsup_{n\to \infty}\dfrac{log(r(n,\epsilon,f))}{n}$. Then \emph{topological entropy of} $f$, $h(f)$ is defined as: \\
\[ h(f)=\displaystyle\lim_{\epsilon\to 0, \epsilon >0} h(\epsilon, f). \]
\end{Definition}

\textbf{Remark:} Note that the Lozi map is defined on $\mathbb R^2$ which is not compact. To be able to investigate the topological entropy of the Lozi map, we take one-point compactification of $\mathbb R^2$ and extend the map continuously to this set. For more details about this continuous extension, see \cite{Ish4}.

\section{Lower Bound Techniques}
There are some computer assisted techniques to give rigorous lower bounds for the topological entropy of maps like H\'enon \cite{Henon} and Ikeda \cite{Ikeda}. They were first introduced by Zygliczy\'{n}ski \cite{Zyg} and developed in \cite{Zyg2} and \cite{Zyg3}. There are also more recent methods by Newhouse, Berz, Makino and Grote \cite{Newhouse} which give better lower bounds for the H\'enon map.

Let us review the following ideas which were used in \cite{Zyg3}.

\begin{figure}[hbtp]
  \vspace{2pt}

  \centerline{\hbox{ 
    \epsfxsize=6.5in
    \epsffile{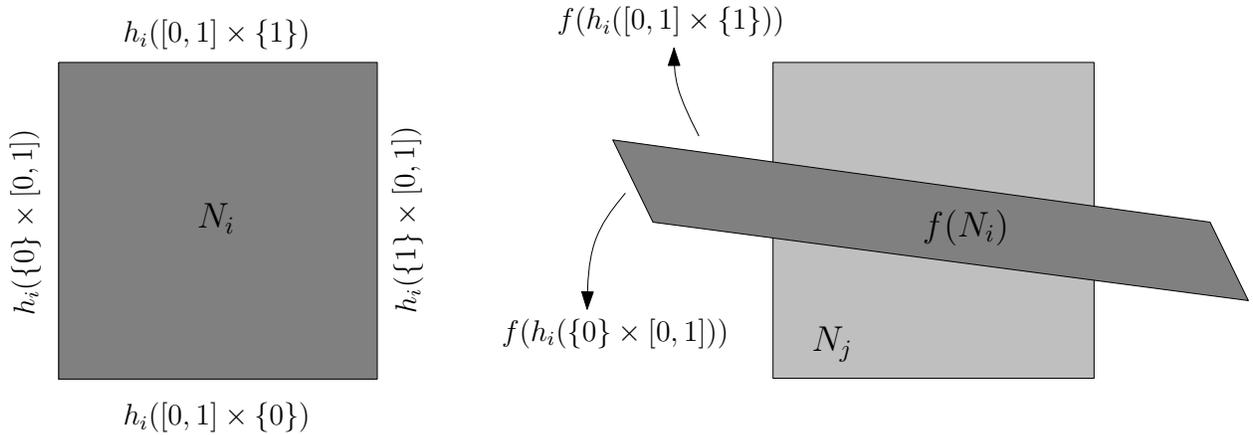}
    }
  }
\caption{$N_i$ $f$-covers $N_j$ \hspace{2mm} ($N_i \Rightarrow N_j$).}
\label{cover}
  \end{figure}

Let $f : \mathbb{R}^2 \to \mathbb{R}^2$ be a continuous map and $N_1, N_2, \dots N_p$ be $p$ pairwise disjoint quadrilaterals. Note that we can parametrize each $N_i$ with the unit square $I^2=[0,1]\times[0,1]$ by choosing a homeomorphism $h_i:I^2 \to N_i$. We call the edges $h_i(\{0\}\times[0,1])$ and $h_i(\{1\}\times[0,1])$ ``vertical" and the edges $h_i([0,1] \times \{0\})$ and $h_i([0,1] \times \{1\})$ ``horizontal". We define a covering relation between two quadrilaterals in the following way (see Fig.~\ref{cover}):

\begin{Definition} We say $N_i$ $f$-covers $N_j$ and write $N_i \Rightarrow N_j$ if:
\begin{itemize}
\item[(i)] $f|N_i$ is one-to-one.
\item[(ii)] For each $\rho \in [0,1]$, there are exactly two numbers $t_{\rho}^1,t_{\rho}^2 \in (0,1)$ such that $f(h_i(\{t_{\rho}^1\}\times\{\rho\}))$ lies in one of the vertical edges of $N_j$ and $f(h_i(\{t_{\rho}^2\}\times\{\rho\}))$ lies in the other vertical edge of $N_j$ and $\forall \thinspace t_{\rho}^1<t<t_{\rho}^2$, $f(h_i(\{t\}\times\{\rho\})) \in N_j$.
\item[(iii)] For $0\leq t < t_{\rho}^1$ and $t_{\rho}^2 < t \leq 1$, $f(h_i(\{t\}\times\{\rho\})) \cap N_j$ is empty.
\end{itemize}
\end{Definition}

If one can show the existence of these quadrilaterals and associated covering relations, they can be used to give rigorous lower bounds for the topological entropy of $f$:
\begin{Theorem}(\cite{Zyg3}) \label{lowerbound} Let $N_1, N_2, \dots N_p$ be pairwise disjoint quadrilaterals and $f : \mathbb{R}^2 \to \mathbb{R}^2$ be continuous. Let $A=(a_{ij})$ be a square matrix where $1\leq i,j \leq p$ and
\[ a_{ij} = \bigg\{ \begin{array}{ccc} 1 & \mbox{ if } & N_i \Rightarrow N_j \\
0 & \mbox{ otherwise }
\end{array}
\]
Then $f$ contains a Cantor set on which it is topologically conjugate to the subshift of finite type with transition matrix $A$. In particular, $h(f)\geq log(\lambda_1)$ where $\lambda_1$ is the largest magnitude eigenvalue ($\lambda_1\geq |\lambda_j|$ for all eigenvalues of $A$).
\end{Theorem}

Note that there is no easy way to detect these systems of quadrilaterals. They are usually found by trial and error. For example, see \cite{Zyg3} and \cite{Newhouse}.

\section{Discontinuity of entropy for Lozi maps}
Buzzi's results \cite{Buzzi} about lower semi-continuity of the entropy of piecewise affine homeomorphisms on compact surfaces can not be applied directly to Lozi maps which are defined on the plane. These results should also hold in the non-compact case, but more work is required. On the other hand, nothing much is known about upper semi-continuity. For Lozi maps, there are some monotonicity results (see \cite{Ish3} and \cite{Burak}) around $b=0$. It is also known that $h(\mathcal L_{a,b})$ depends continuously on the parameters $(a,b)$ at all points $(a,0)$ where $a>1$: First note that $h(\mathcal L_{a,0})= min\{loga,log2\}$ for $a>1$ as in the tent map. By the monotonicity results in \cite{Ish3}, $h(\mathcal L_{a-N|b|,0})\leq h(\mathcal L_{a,b})\leq h(\mathcal L_{a+N|b|,0})$ for some $N>0$ and $|b|$ small. So continuity follows. \\

We first prove that the entropy jumps from zero to a positive value if parameters are slightly changed from $(a,b)=(1.4, 0.4)$ to $(a,b)=(1.4+\epsilon_1, 0.4)$ where $\epsilon_1$ is positive and small.

\begin{Theorem}\label{Preliminary A} There exists some $\epsilon_*>0$ such that for all $0<\epsilon_1<\epsilon_*$:
\begin{itemize}
\item[(i)] The topological entropy of the Lozi map with $(a,b)=(1.4,0.4)$, $h(\mathcal L_{1.4,0.4})$, is zero.
\item[(ii)] The topological entropies of the Lozi maps, $h(\mathcal L_{(1.4+\epsilon_1,0.4)})$, have a lower bound of $0.1203$.
\end{itemize}

\end{Theorem}

\begin{proof}[Proof of Theorem \ref{Preliminary A} (i)] \emph{ }\\
Let us denote $\mathcal L_{1.4, 0.4} = \mathcal L$. We will prove that $h(\mathcal L^4)=0$.
By direct calculation of $\mathcal L^4$, one can solve the equation $\mathcal L^4(x,y)=(x,y)$ for $(x,y) \in \mathbb{R}^2$ to see that $\mathcal L^4$ has the following fixed points (see the Appendix):
\begin{itemize}
\item[$(i)$] Fixed points of $\mathcal L$: $p_1=(1/2, 1/2)$ and $p_2=(-5/4, -5/4)$,
\item[$(ii)$] The closed line segment $\ell_1$ which connects $(-20/29, 35/29)$ to \\
$(0,15/29)=\mathcal L^2(-20/29,35/29)$,
\item[$(iii)$] The closed line segment which connects $(15/29,-20/29)$ to \\
$(35/29,0)=\mathcal L^2(15/29,-20/29)$, i.e.~$\mathcal L(\ell_1)$.
\end{itemize}

\noindent Note that $p_1$ is a saddle fixed point and $v^{s}_1=(\lambda^{s}_1, 1)$ where $\lambda^{s}_1=(-7+\sqrt{89})/10$ is a stable direction at $p_1$ and $W^{s}_{+}(p_1)=\{p_1+v^{s}_{1}t \in \mathbb{R}^2 | \thinspace t>0\}$ is invariant under $\mathcal L$ (and therefore $\mathcal L^4$). Similarly, $p_2$ is a saddle point and $v^{u}_2=(-\lambda^{u}_2, -1)$ where $\lambda^{u}_2=(7+\sqrt{89})/10$ is an unstable direction at $p_2$ and $W^{u}_{+}(p_2)=\{p_2+v^{u}_{2}t \in \mathbb{R}^2 | \thinspace t>0\}$ is invariant under $\mathcal L^4$. \\

\noindent Let us call the left and the right connected components of the unstable manifold at $p_1$; $W_{\ell}(p_1)$ and $W_r(p_1)$, respectively (see Fig.~\ref{Upicture}). We want to show that $W_{\ell}(p_1)$ is attracted by $\ell_1$ and $W_r(p_1)$ is attracted by $\mathcal L(\ell_1)$. But let us first explain how to conclude the proof of Theorem \ref{Preliminary A} $(i)$ from that claim. Let $U = \mathbb{R}^2\setminus M$ where $M=W^s_{+}(p_1)\cup\{p_1\}\cup W^u_{+}(p_2)\cup\{p_2\}\cup W_\ell(p_1)\cup\ell_1\cup W_r(p_1)\cup\mathcal L(\ell_1)$. Note that $U$ is invariant by construction and it is simply connected since the complement of $U$ in the extended plane, i.e.~$M\cup \{\infty\}$, is connected. Also, note that $M\cup \{\infty\}$ is compact because of the claim that $W_{\ell}(p_1)$ is attracted by $\ell_1$ and $W_r(p_1)$ is attracted by $\mathcal L(\ell_1)$. This implies $U$ is homeomorphic to the open unit disk (by Riemann Mapping Theorem) which is homeomorphic to $\mathbb{R}^2$ .  Since $\mathcal L^4$ has no fixed points in $U$ and it is orientation preserving, Brouwer's translation theorem implies that $\mathcal L^4$ has no non-wandering points in $U$. This shows the non-wandering set of $\mathcal L^4$ only consists of the fixed points of $\mathcal L^4$. So, $h(\mathcal L^4)=4h(\mathcal L)=0$.

\begin{figure}[hbtp]
  \centerline{\hbox{ 
    \epsfxsize=11in
    \epsfysize=0in
    \epsffile{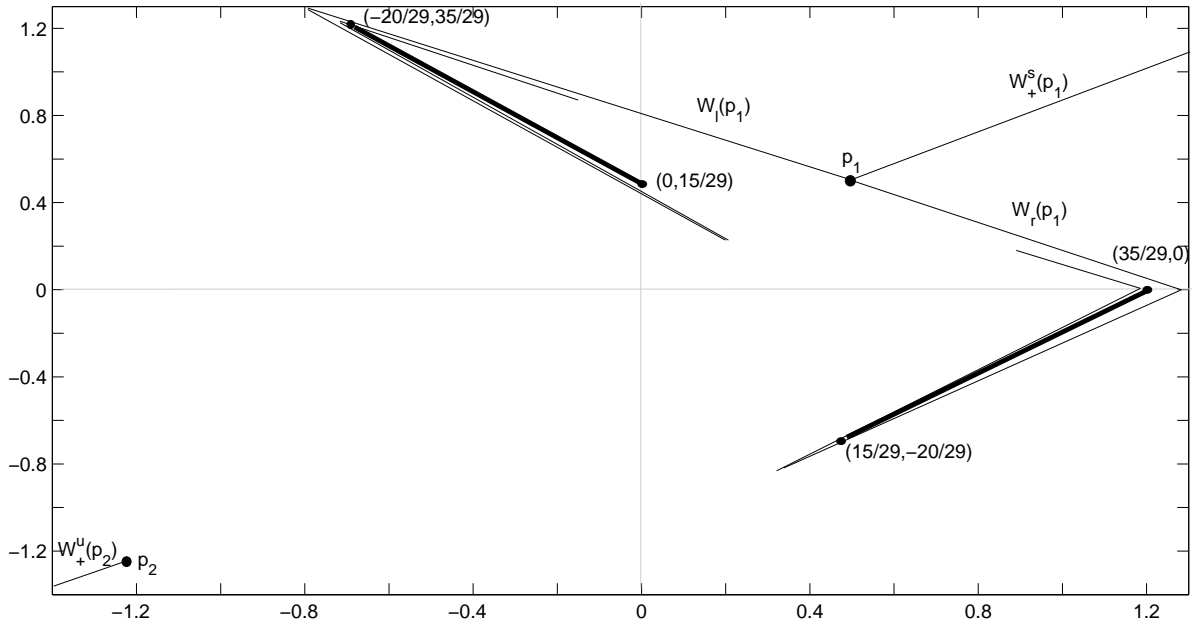}
    }
  }
\caption{Several components of the invariant manifolds of the fixed points $p_1$ and $p_2$ given together with the two line segments (darker) of period-$4$ points: $\ell_1$ which connects $(-20/29, 35/29)$ to $(0,15/29)$ and $\mathcal L(\ell_1)$ which connects $(15/29,-20/29)$ to $(35/29,0)$. $U$ is the complement of all the points shown in the picture.}
\label{Upicture}
  \end{figure}

%
\begin{figure}[hbtp]
  \centerline{\hbox{ 
    \epsfxsize=9in
    \epsffile{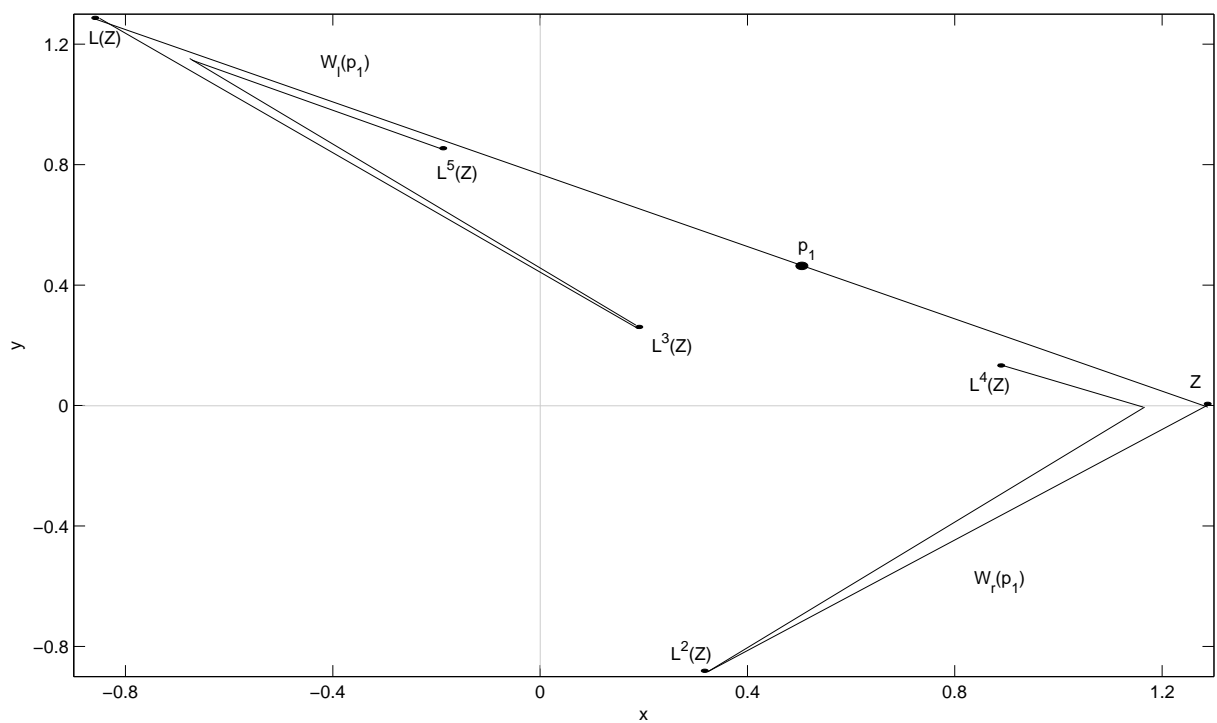}
    }
  }
\caption{The left and right connected components of the unstable manifold of $p_1$ and several iterations of the point $Z$ under the map $\mathcal L = \mathcal L_{1.4,0.4}$.}
\label{fullpicture}
  \end{figure}

\paragraph{$W_{\ell}(p_1)$ is attracted to $\ell_1$:} Now, let $Z$ be the intersection of the half line $m = \{p_1+v_{1}^ut \in \mathbb{R}^2 \thinspace | \thinspace t>0 \}$ and the $x$-axis where $v_{1}^u=(-\lambda_{1}^u, -1)$ and $\lambda_1^u=(-7-\sqrt{89})/10$ (see Fig.~\ref{fullpicture}). Note that $W_{\ell}(p_1)=\bigcup_{n=0}^{\infty}\mathcal L^{4n}(\{p_1-v_{1}^ut \thinspace | \thinspace 0.1>t>0\})$, i.e.~forward iterations of a small piece in the unstable direction. Let the portion of $W_{\ell}(p_1)$ which connects $\mathcal L(Z)$ and $\mathcal L^5(Z)$ be called $W$. It is not hard to see that $W_{\ell}(p_1)=\bigcup_{n=-\infty}^{\infty}\mathcal L^{4n}(W)$. We want to show that every $x \in W$ (so every $x \in W_{\ell}(p_1)$) is attracted to $\ell_1$.\\

\textbf{Remark:} Note that all points in $\ell_1$ have a neutral direction (along $\ell_1$) and a contracting
direction with slope $-5/2$. This gives an immediate basin of attraction up to the interaction with the
singularity lines of $\mathcal L^4$. The basin (trapping region) intersects and therefore captures a large part of $W_{\ell}(p_1)$
but not all since that set extends to the left and right. Below, we show that these left and right parts are also eventually attracted to the trapping region. \\

\begin{figure}[hbtp]

  \centerline{\hbox{ 
    \epsfxsize=6.5in
    \epsffile{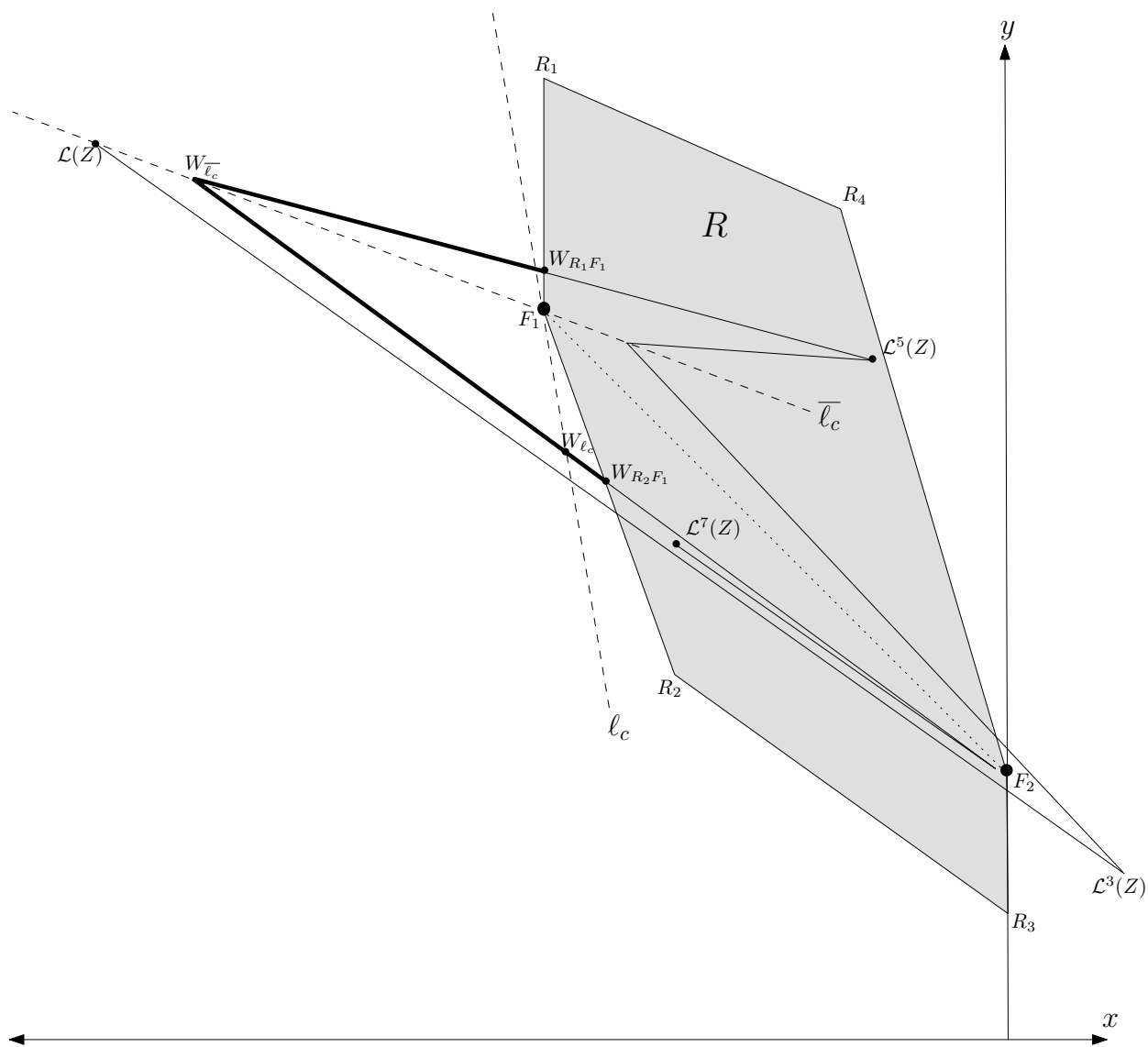}
    }
  }
\caption{This figure shows a portion of the left unstable manifold of the fixed point $p_1$. Note that all the points on the line segment connecting
$F_1$ to $F_2$ are period-$4$ points of $\mathcal L$ }
\label{left_unstable}
  \end{figure}

\noindent \emph{Trapping Region:} Let $f:\mathbb{R}^n \to \mathbb{R}^n $ be a map. A neighborhood $U$ of an $f$-invariant set $A\subset\mathbb{R}^n$ is called a \emph{trapping region for A} if $\forall m >0$, $f^m(U) \subset U$ and $\bigcap_{m>0}f^m(U)=A$. Below, we introduce a trapping region $R$ around $\ell_1$ such that any point $x\in R$ is attracted to a point in $\ell_1$ under forward iterations of $\mathcal L$. Let:
\begin{align*}
& R_1=(-20/29, 35/29 + 0.2) \\
& R_2=(-20/29 + 0.1, 35/29 - 0.25) \\
& R_3=(0, 15/29 - 0.25) \\
& R_4=(-0.2, 15/29 + 0.5)
\end{align*}
\noindent Let us call the left and right end points of $\ell_1$; $F_1$ and $F_2$, respectively. Note that $F_1=(-20/29,35/29)$ and $F_2=\mathcal L^2(F_1)=(0,15/29)$. Let $R$ be the hexagon with vertices $R_1$,$F_1$,$R_2$,$R_3$,$F_2$ and $R_4$. The sides $F_1R_2$ and $F_2R_4$ are parallel to each other with slope $-5/2$ and they are stable directions at $F_1$ and $F_2$, respectively. Since $R_1$ is in the stable manifold of a point in $\ell_1$, it is attracted to $\ell_1$ under iterations of $\mathcal L^4$. Similarly, $R_4$ is attracted to $F_2$ since it is in the stable manifold of $F_2$. So, the quadrilateral with vertices $R_1$,$F_1$,$F_2$ and $R_4$ is mapped to thinner and thinner quadrilaterals for which one of the sides is always $\ell_1=F_1F_2$. Similarly, the quadrilateral with vertices $F_1$,$R_2$,$R_3$ and $F_2$ is mapped towards $\ell_1$ (see Fig.~\ref{trapping}). So, $R$ is a trapping region. \\

\begin{figure}[hbtp]

  \centerline{\hbox{ 
    \epsfxsize=3.5in
    \epsffile{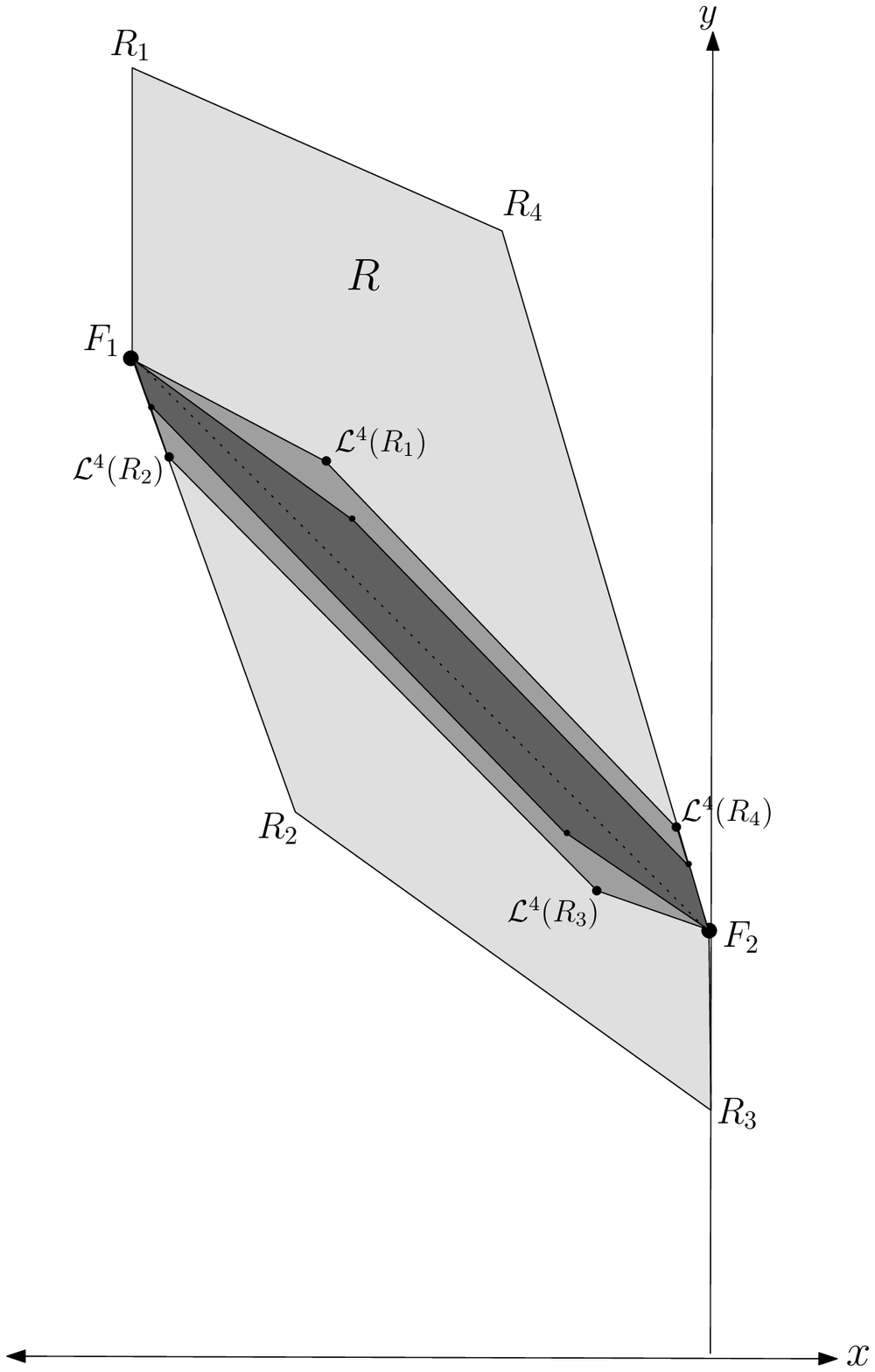}
    }
  }
\caption{Trapping region $R$ (gray) and images $\mathcal L^4(R)$ (darker) and $\mathcal L^8(R)$ (darkest). }
  \label{trapping}
  \end{figure}

\noindent We want to show that all the points in $W$ are eventually mapped into $R$ under forward iterations of $\mathcal L^4$.
Let us start with the part of $W$ which connects $\mathcal L(Z)$ and $\mathcal L^3(Z)$.
The image of this line segment (under $\mathcal L^4$) is the portion of $W_{\ell}(p_1)$ which connects $\mathcal L^5(Z)$ and
$\mathcal L^7(Z)$ (see Fig.~\ref{left_unstable}).
Let us call this portion $\overline{W}$. $\mathcal L^5(Z)$ and $\mathcal L^7(Z)$ are both in $R$ but there is a part of $\overline{W}$ which is still
outside of $R$ which we denote by $\overline{\overline{W}}$, i.e.~$\overline{\overline{W}}$ is the closure of $\overline{W}\setminus R$. Note that $\ell_c : y=1-1.4(1+1.4x+0.4y)+0.4x$ is a critical line for $\mathcal L^4$ around $F_1$, i.e. images of lines which transversally intersect $\ell_c$  are broken lines. Let $\overline{\ell_c}=\mathcal L^4(\ell_c)$. Also, let $W\cap R_1F_1=W_{R_1F_1}$,
$W\cap R_2F_1= W_{R_2F_1}$, $W\cap \overline{\ell_c}=W_{\overline{\ell_c}}$ and the intersection point of $W$ and $\ell_c$ which stays below
$\overline{\ell_c}$ be $W_{\ell_c}$. \\

\begin{figure}[hbtp]

  \centerline{\hbox{ 
    \epsfxsize=5.5in
    \epsffile{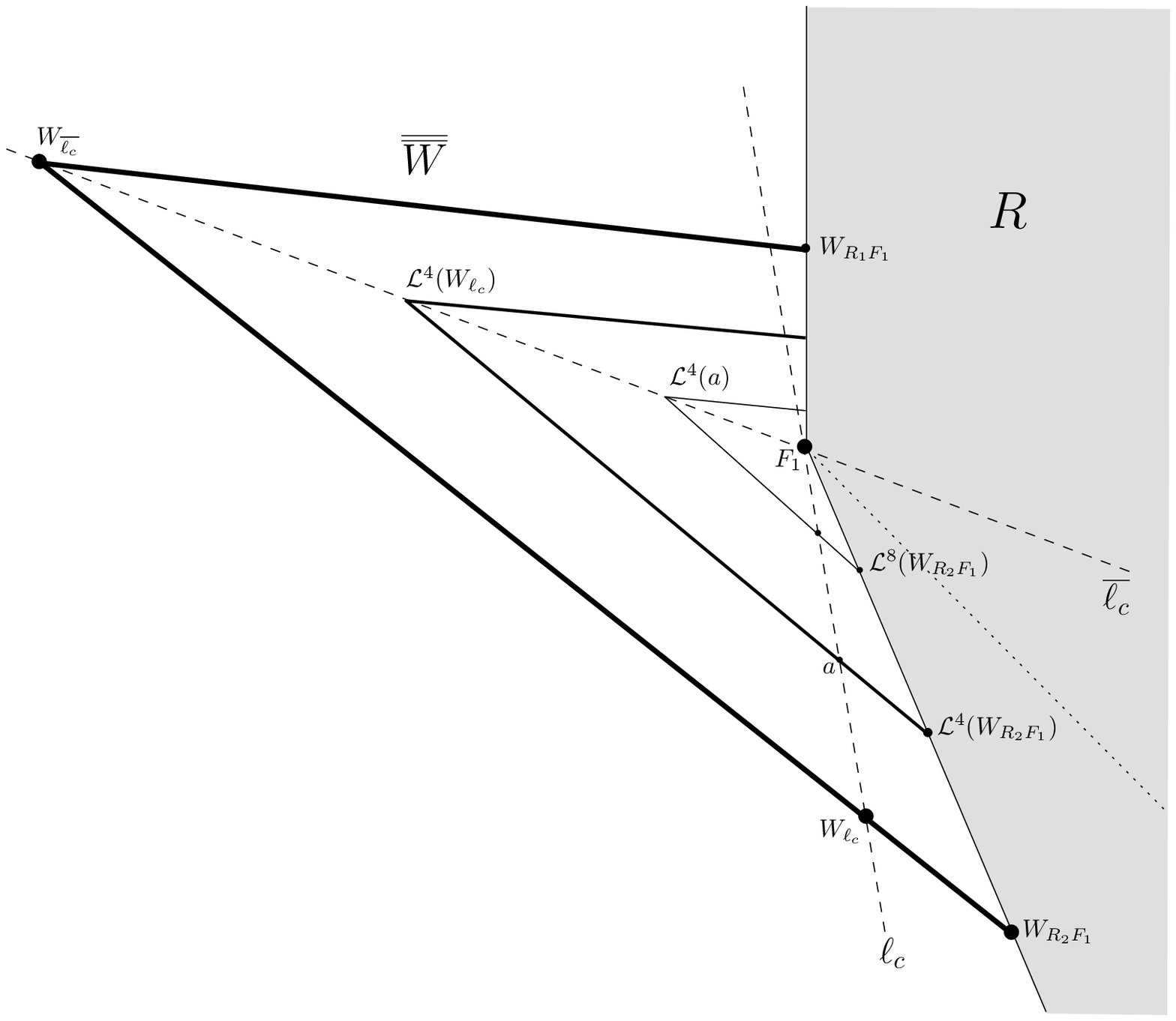}
    }
  }
\caption{The set $\overline{\overline{W}}$ (thickest solid broken line) and the part of the images $\mathcal L^4(\overline{\overline{W}})$ (thinner) and $\mathcal L^8(\overline{\overline{W}})$ (thinnest) which stay outside of $R$. Note that everything above $\overline{\ell_c}$ is mapped into $R$ under $\mathcal L^4$.}
  \label{left_unstable_zoom}
  \end{figure}

\noindent $\overline{\overline{W}}$ consists of two parts: The line segment which connects $W_{R_1F_1}$ and $W_{\overline{\ell_c}}$ and the line segment which connects $W_{\overline{\ell_c}}$ and $W_{R_2F_1}$ (see Fig.~\ref{left_unstable_zoom}). Note that $\mathcal L^4(\overline{\ell_c})$ is a broken line that stays in $R$ since $\overline{\ell_c}$ intersects $\ell_c$ which is a critical line for $\mathcal L^4$. So, all points on the line segment connecting $W_{R_1F_1}$ and $W_{\overline{\ell_c}}$ are mapped into $R$, too. \\
On the other hand, $W_{\ell_c}$ is mapped to a point on $\overline{\ell_c}$. So, the line segment connecting $W_{\ell_c}$ and $W_{\overline{\ell_c}}$ is also completely mapped into $R$ under $\mathcal L^8$.\\
The only part left is the portion that connects $W_{\ell_c}$ and $W_{R_2F_1}$. But note that $W_{R_2F_1}$ is on the stable direction so forward iterations move towards $F_1$. $W_{\ell_c}$ is mapped between $W_{\overline{\ell_c}}$ and $F_1$. So, one can repeat the same argument to this line segment connecting $\mathcal L^4(W_{R_2F_1})$ and $\mathcal L^4(W_{\ell_c})$. So, by induction the portion that connects $W_{\ell_c}$ and $W_{R_2F_1}$ is also mapped into $R$ eventually. This completes the proof that $\overline{\overline{W}}$ is mapped into $R$. \\
The above analysis explains that forward images of $\overline{\overline{W}}$ consists of some parts which is mapped into $R$ and some parts which stays outside of $R$. However, the parts outside of $R$ are eventually attracted by $R$ (see Fig.~\ref{left_unstable_zoom}). \\

\noindent Now, for the other portion of $W$ (connecting $\mathcal L^3(Z)$ and $\mathcal L^5(Z)$) similar argument can be applied while this time the critical line $\ell_c$ is the $y$-axis and the parts outside of $R$ are either mapped into $R$ or attracted by $F_2$.\\

\noindent Finally, note that $W_{\ell}(p_1)$ is attracted to $\ell_1$ implies that $W_{r}(p_1)=\mathcal L(W_{\ell}(p_1))$ is attracted to $\mathcal L(\ell_1)$.
\end{proof}

\begin{proof}[Proof of Theorem \ref{Preliminary A} (ii)] \emph{ } \\
We want to show that for any $\epsilon_1$ positive and small, Theorem \ref{lowerbound} applies with an appropriate
subshift of finite type yielding the lower bound for the map $\mathcal L_{(1.4+\epsilon_1,0.4)}$.

\noindent Fix an $\epsilon_1>0$ and denote $\mathcal L_{\epsilon_1}=\mathcal L_{(1.4+\epsilon_1,0.4)}$. Note that the line segment
connecting $F_1=(-20/29,35/29)$ and $F_2=(0,15/29)$ consists of period-$4$ points of $\mathcal L_{(1.4,0.4)}$.\\

\noindent Now, let $N_1$ be the quadrilateral given by the four vertices:
\begin{align*}
& A=(0, 15/29 - \epsilon_1) \\
& B=(\epsilon_1, 15/29 + (7/2)\epsilon_1) \\
& C=((5/2)\epsilon_1, 15/29 + (5/2)\epsilon_1) \\
& D=((3/2)\epsilon_1, 15/29 - 2\epsilon_1)
\end{align*}

Also let $N_2$ be the quadrilateral whose vertices are:
\begin{align*}
& E=(-3\epsilon_1, 15/29 + (7/2)\epsilon_1) \\
& F=(-2\epsilon_1, 15/29 + (5/6)\epsilon_1) \\
& G=(0, 15/29 - (1/2)\epsilon_1) \\
& H=(-\epsilon_1, 15/29 + (13/6)\epsilon_1)
\end{align*}

\noindent For $N_1$, let the sides $AB$ and $CD$ be ``vertical" and the other two sides be ``horizontal". Similarly for $N_2$, let $EF$ and $GH$ be ``vertical" and the other two sides be ``horizontal". Note that the images of $N_1$ and $N_2$ under $\mathcal L_{\epsilon_1}^4$ are also quadrilaterals since $N_1$ and $N_2$ are chosen away from the singularity locus of $\mathcal L_{\epsilon_1}^4$. Moreover, vertical edges are contracted since they are close to the stable directions around $(0,15/29)$ and $(-20/29,35/29)$.\\

\noindent By direct calculation, it can be shown that the images of the vertices under the map $\mathcal L_{\epsilon_1}^4$ is given by (see Fig.~\ref{boxes}):
\begin{align*}
& \mathcal L_{\epsilon_1}^4(A)=(\frac{30476}{18125}\epsilon_1 + O(\epsilon_1 ^2), \frac{15}{29} - \frac{6363}{3625}\epsilon_1 + O(\epsilon_1^2))
\approx (1.68\epsilon_1, \frac{15}{29}-1.75\epsilon_1) \\
& \mathcal L_{\epsilon_1}^4(B)=(\frac{6188}{3625}\epsilon_1 + O(\epsilon_1 ^2), \frac{15}{29} - \frac{1319}{725}\epsilon_1 + O(\epsilon_1^2))
\approx (1.70\epsilon_1, \frac{15}{29}-1.81\epsilon_1) \\
& \mathcal L_{\epsilon_1}^4(C)=(-\frac{4769}{1450}\epsilon_1 + O(\epsilon_1 ^2), \frac{15}{29} + \frac{847}{290}\epsilon_1 + O(\epsilon_1^2))
\approx (-3.28\epsilon_1, \frac{15}{29}+2.92\epsilon_1) \\
& \mathcal L_{\epsilon_1}^4(D)=(-\frac{120153}{36250}\epsilon_1 + O(\epsilon_1 ^2), \frac{15}{29} + \frac{21639}{7250}\epsilon_1 + O(\epsilon_1^2))
\approx (-3.31\epsilon_1, \frac{15}{29}+2.98\epsilon_1) \\
& \mathcal L_{\epsilon_1}^4(E)=(-\frac{9283}{18125}\epsilon_1 + O(\epsilon_1 ^2), \frac{15}{29} + \frac{1554}{3625}\epsilon_1 + O(\epsilon_1^2))
\approx (-0.51\epsilon_1, \frac{15}{29}+0.42\epsilon_1) \\
& \mathcal L_{\epsilon_1}^4(F)=(-\frac{23209}{54375}\epsilon_1 + O(\epsilon_1 ^2), \frac{15}{29} + \frac{3792}{10875}\epsilon_1 + O(\epsilon_1^2))
\approx (-0.42\epsilon_1, \frac{15}{29}+0.34\epsilon_1) \\
& \mathcal L_{\epsilon_1}^4(G)=(\frac{36363}{18125}\epsilon_1 + O(\epsilon_1 ^2), \frac{15}{29} - \frac{7494}{3625}\epsilon_1 + O(\epsilon_1^2))
\approx (2.00\epsilon_1, \frac{15}{29}-2.06\epsilon_1) \\
& \mathcal L_{\epsilon_1}^4(H)=(\frac{113584}{54375}\epsilon_1 + O(\epsilon_1 ^2), \frac{15}{29} - \frac{22917}{10875}\epsilon_1 + O(\epsilon_1^2))
\approx (2.08\epsilon_1, \frac{15}{29}-2.10\epsilon_1) \\
\end{align*}

\begin{figure}[htbp]

  \centerline{\hbox{ 
    \epsfxsize=5.5in
    \epsffile{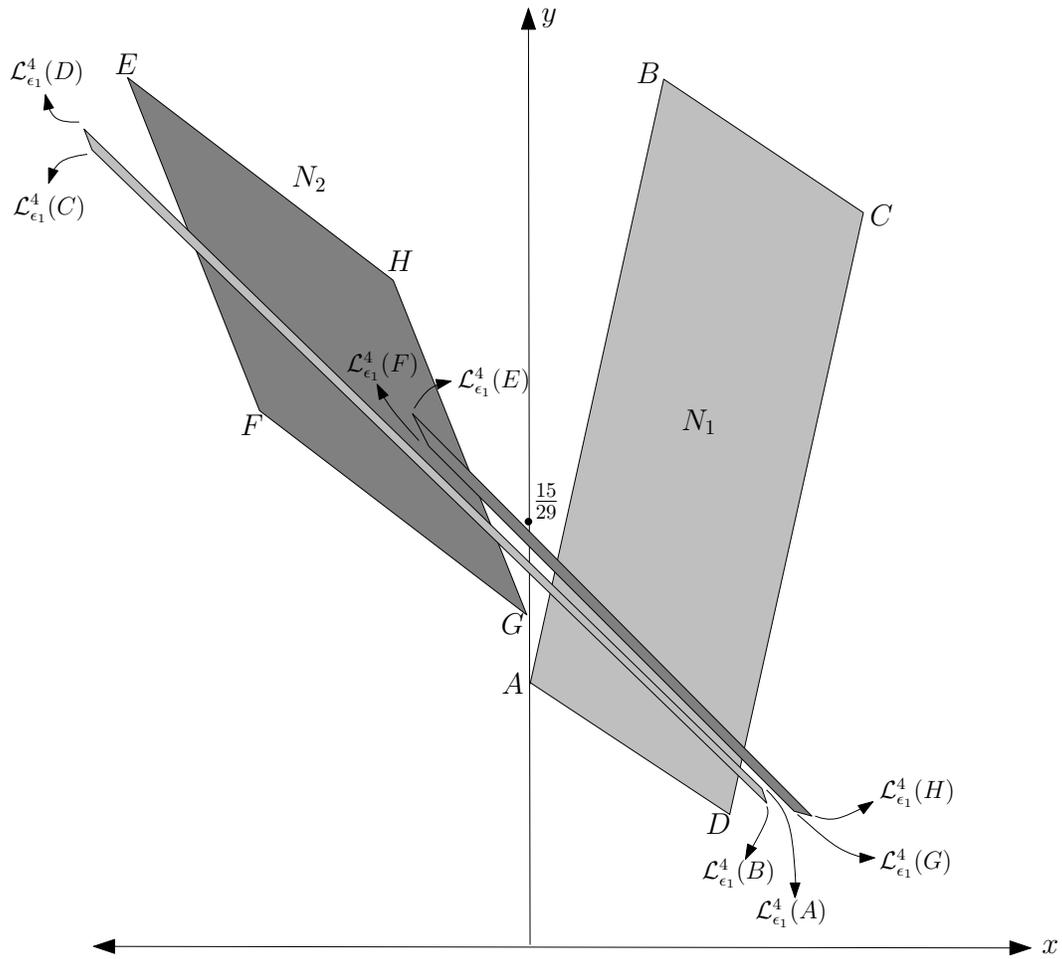}
    }
  }
\caption{This figure shows the quadrangles $N_1$ and $N_2$ and their images (thinner boxes). Notice the covering relations:
  $N_1 \Rightarrow N_1$, $N_1 \Rightarrow N_2$ and $N_2 \Rightarrow N_1$ }
  \label{boxes}
  \end{figure}
It is not hard to see that we have the following covering relations: $N_1 \Rightarrow N_1$, $N_1 \Rightarrow N_2$ and $N_2 \Rightarrow N_1$.
So the transition matrix is given by:
\[ \left( \begin{array}{ccc}
1 & 1 \\
1 & 0 \end{array} \right)
\]
where the largest magnitude eigenvalue is $\frac{\sqrt{5}+1}{2}$. Since we are using $\mathcal L_{\epsilon_1}^4$ during the process $h(\mathcal L_{\epsilon_1})$=
$\frac{1}{4}h(\mathcal L_{\epsilon_1}^4)\geq \frac{1}{4} log\frac{\sqrt{5}+1}{2} > 0.1203$ by Theorem \ref{lowerbound}.
\end{proof}

\begin{figure}[htbp]
   \centerline{\hbox{
    \epsfxsize=4.5in
    \epsfysize=3in
    \epsffile{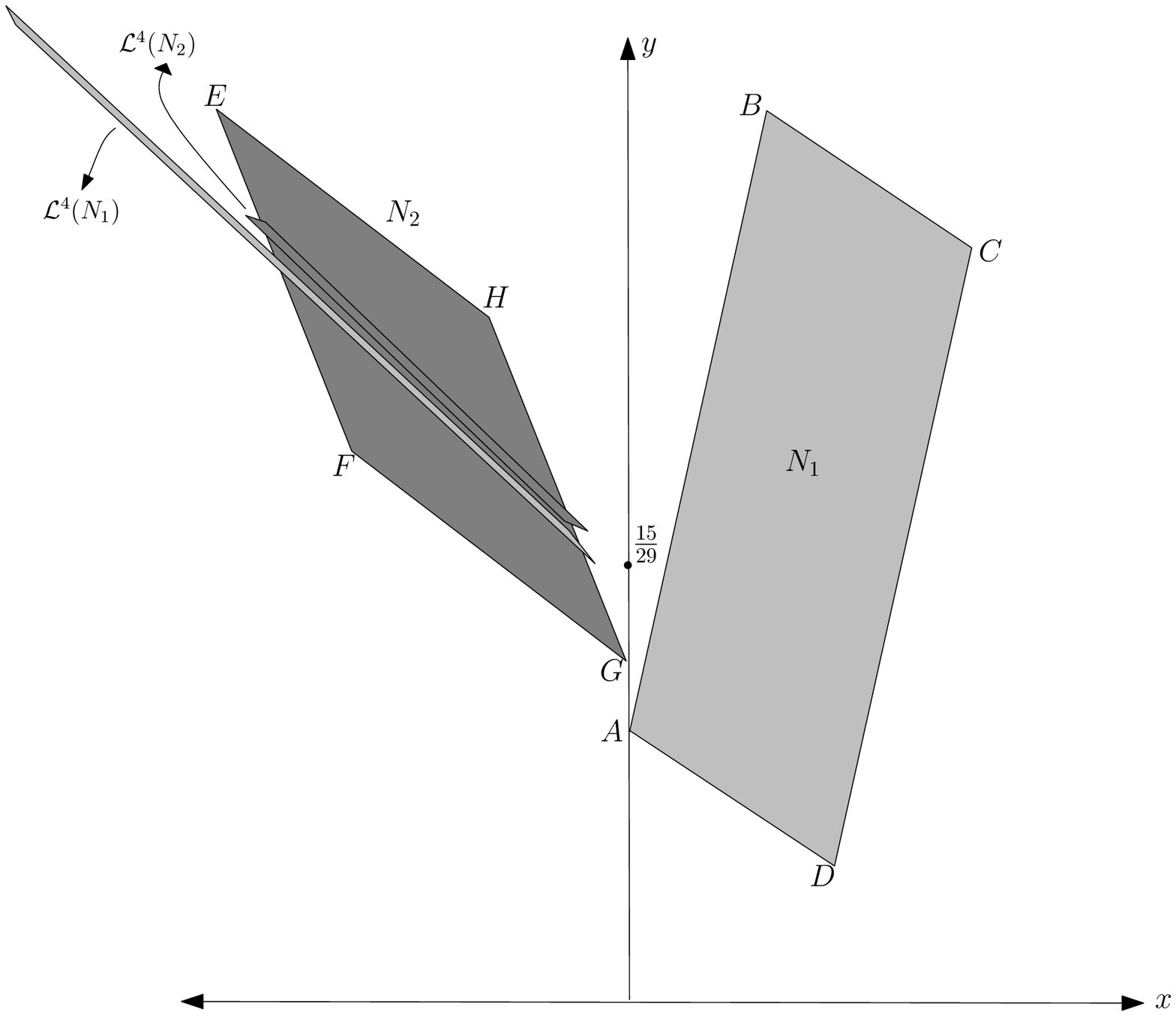}
    \hspace{0.55in}
    }
    }
    \vspace{1in}
    \centerline{ \hbox{
    \epsfxsize=4in
     \epsfysize=3in
    \epsffile{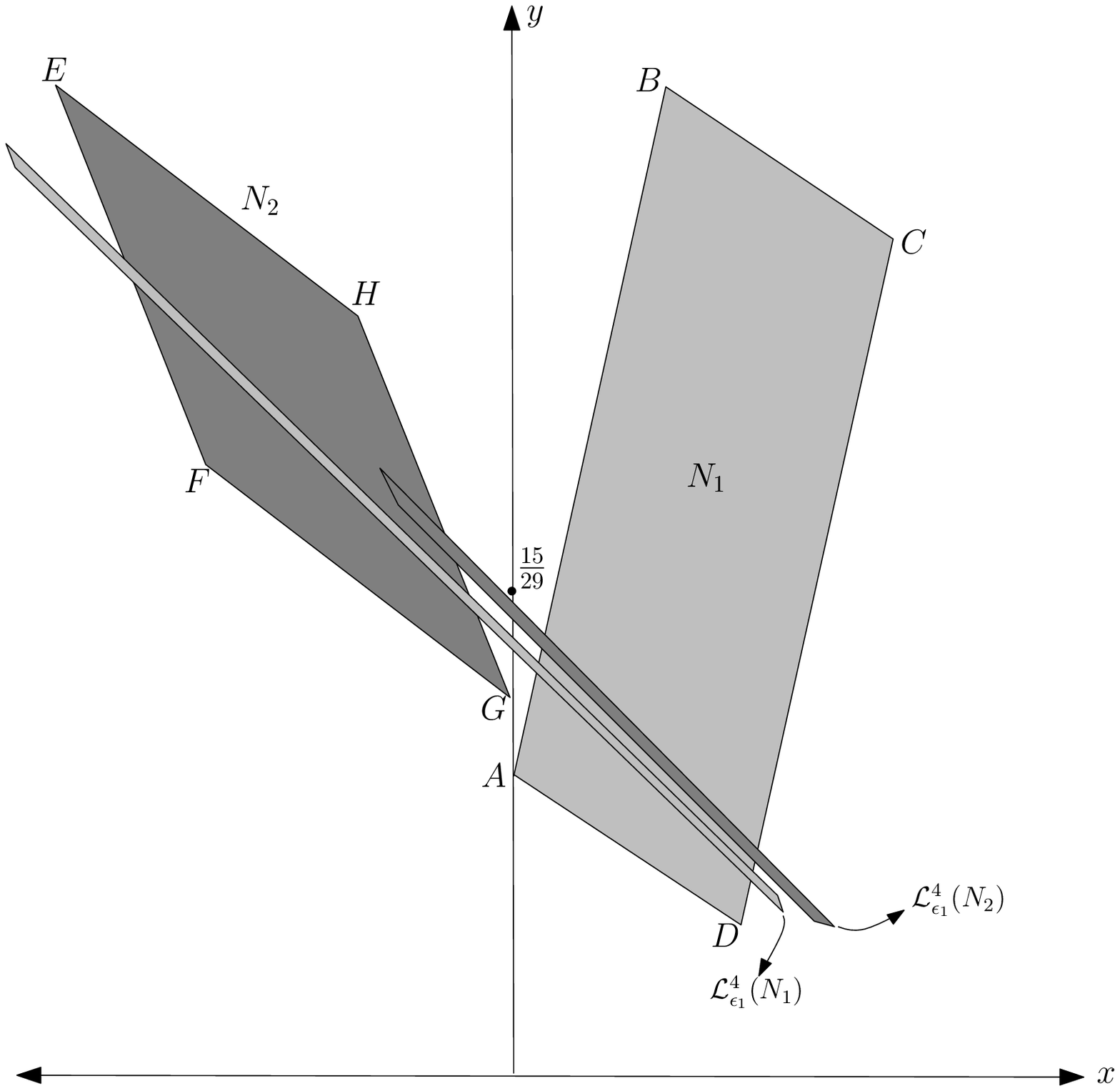}
    }
  }
  \parbox{6in}{ \hspace{0.2in}
  \caption{The comparison between the images of $N_1$ and $N_2$ under the maps $\mathcal L^4$ (top) and $\mathcal L_{\epsilon_1}^4$ (bottom). $N_1$ and $N_2$ are the same in both cases with some small, fixed $\epsilon_1>0$. Note that both pictures are in a small neighborhood of $(0,15/29)$. There is not enough expansion in the first case (with $\mathcal L^4$)  to create covering relations but the perturbed map $\mathcal L_{\epsilon_1}^4$ creates enough expansion which causes a jump up at the entropy. \label{comparison}} }

\end{figure}

\textbf{Remark:} We would like to point out that the jump up in the entropy explained above is somewhat similar to the following one dimensional case: Let $T:\mathbb{R}\to\mathbb{R}$ be defined by $T(x)=-2|x|$. All the initial points except the fixed point at $x=0$ go to infinity under further iterations of $T$ so the entropy of $T$ is zero. Note that the graph of $T(x)$ stays below the diagonal line $y=x$. On the other hand, the perturbed map $T_{\delta}(x)=-2|x|+\delta$ where $\delta>0$ has entropy $log2$ (similar to the standard tent map) and the graph of $T_{\delta}(x)$ crosses the diagonal line. One can see a similar kind of behavior at the images of $N_1$ and $N_2$ under the maps $\mathcal L^4$ and $\mathcal L_{\epsilon_1}^4$ (see Fig.~\ref{comparison}). Images of $N_1$ and $N_2$ under $\mathcal L^4$ stay on the left of the critical line $x=0$ and the entropy is zero. On the other hand, under $\mathcal L_{\epsilon_1}^4$, these images cross the critical line and the entropy jumps up. We would like to thank S.~E.~Newhouse for pointing out this similarity between the one dimensional and two dimensional cases. \\

Now, we can extend our results from $(a,b)=(1.4,0.4)$ to $(a,b)=(1.4+\epsilon_2, 0.4+\epsilon_2)$ where $|\epsilon_2|$ is small: \\

\begin{proof}[{\bf Proof of Theorem \ref{A} }] \emph{ } \\
Let $\mathcal L$ denote $\mathcal L_{(1.4+\epsilon_2, 0.4+\epsilon_2)}$. \\

\emph{(i) The entropy is zero for $\mathcal L$:} \\

For $|\epsilon_2|$ small and fixed, we still have two line segments of period-4 points: the line segment connecting $F_2^{\epsilon_2}=\frac{1-(0.4+\epsilon_2)^2}{(1.4+\epsilon_2)(1+(0.4+\epsilon_2)^2)}$ and $F_1^{\epsilon_2}=\mathcal L^2(F_2^{\epsilon_2})$ and the image of this line segment under $\mathcal L$. So, we can still find a similar trapping region using the vertical lines and the stable directions at $F_1^{\epsilon_2}$ and $F_2^{\epsilon_2}$. The rest of the proof is the same as in the case of $(a,b)=(1.4,0.4)$. \\

\emph{(ii) The lower bound for $(a,b)=(1.4+\epsilon_1+\epsilon_2,0.4+\epsilon_2)$:}

Let $\mathcal L_{\epsilon_1}=\mathcal L_{(1.4+\epsilon_1+\epsilon_2,0.4+\epsilon_2)}$. We need to find two boxes as in the case of $(a,b)=(1.4,0.4)$ which give us the covering relations. We slightly modify the points we used before: \\
For $\epsilon_1$ positive and small, let $N_1$ be the quadrilateral given by the four vertices:
\begin{align*}
& \tA=(0, F_2^{\epsilon_2} - \epsilon_1) \\
& \tB=(\epsilon_1, F_2^{\epsilon_2} + (7/2)\epsilon_1) \\
& \tC=((5/2)\epsilon_1, F_2^{\epsilon_2} + (5/2)\epsilon_1) \\
& \tD=((3/2)\epsilon_1, F_2^{\epsilon_2} - 2\epsilon_1)
\end{align*}

Also let $N_2$ be the quadrilateral whose vertices are:
\begin{align*}
& \tE=(-3\epsilon_1, F_2^{\epsilon_2} + (7/2)\epsilon_1) \\
& \tF=(-2\epsilon_1, F_2^{\epsilon_2} + (5/6)\epsilon_1) \\
& \tG=(0, F_2^{\epsilon_2} - (1/2)\epsilon_1) \\
& \tH=(-\epsilon_1, F_2^{\epsilon_2} + (13/6)\epsilon_1)
\end{align*}

In other words, $15/29$ is replaced with $F_2^{\epsilon_2}$. We want to show that we still have the same covering relations and the same lower bound.  \\
Although one can explicitly write down the images of these points under $\mathcal L_{\epsilon_1}^4$, for simplicity we only want to point out the differences between this case and the case $(a,b)=(1.4,0.4)$. For example, $\mathcal L_{\epsilon_1}^4(\tA)$ consists of terms including $\epsilon_1$ and some others not including $\epsilon_1$. Observe that if $\epsilon_1$ equals zero then $F_2^{\epsilon_2}$ is a period-$4$ point, so the terms in $\mathcal L_{\epsilon_1}^4(\tA)$ \emph{not including} $\epsilon_1$ add up to $F_2^{\epsilon_2}$ (This is because when $\epsilon_1 = 0$, $\tA$ becomes $(0, F_2^{\epsilon_2})$ and so $\mathcal L_{(1.4 + \epsilon_2, 0.4+\epsilon_2)}^4(\tA) = F_2^{\epsilon_2}$). Note that $15/29$ in the proof of $(a,b)=(1.4,0.4)$ case is now replaced with $F_2^{\epsilon_2}$.  \\
On the other hand, the terms in $\mathcal L_{\epsilon_1}^4(\tA)$ \emph{including} $\epsilon_1$ can be made arbitrarily close to the terms including $\epsilon_1$ in the $(a,b)=(1.4,0.4)$ case (i.e.~to the terms $(30476/18125)\epsilon_1$ in the $x$-coordinate and $-(6363/3625)\epsilon_1$ in the $y$-coordinate of $\mathcal L_{\epsilon_1}^4(A)$) by choosing $|\epsilon_2|$ small. Note that the size of $|\epsilon_2|$ does not depend on $\epsilon_1$ but rather depends on the coefficient of $\epsilon_1$, i.e.~$(30476/18125)$ and $-(6363/3625)$. \\
The same argument can be applied to all other points, so our new boxes also satisfy the previous covering relations giving the same lower bound ($0.1203$) for the entropy.

\end{proof}

\textbf{Remark:} The reason why the entropy is zero on the line segment $l=\{(1.4 + \epsilon_2, 0.4 +\epsilon_2): |\epsilon_2|<\epsilon_*\}$ and it is above $0.1203$ for the parameters to the right of that segment is the fact that we have a line segment of period-$4$ points when the parameters are chosen from $l$. These period-$4$ points create a trapping region causing the zero entropy. On the other hand, period-$4$ points suddenly disappear to the right of $l$ causing enough expansion and allowing us to find the necessary subshift which gives the positive entropy (see Fig.~\ref{comparison}).

\appendix
\section{Appendix}
Here, we explain some of the details in the proof of Theorem \ref{Preliminary A} $(i)$. We show that $\mathcal L_{1.4,0.4}^4 = \mathcal L^4$ has the following fixed points: $(i)$ fixed points of $\mathcal L_{1.4,0.4}=\mathcal L$: $p_1=(1/2, 1/2)$ and $p_2=(-5/4, -5/4)$, $(ii)$ the closed line segment $\ell_1$ which connects $(-20/29, 35/29)$ to $(0,15/29)=\mathcal L^2(-20/29,35/29)$ and $(iii)$ the closed line segment which connects $(15/29,-20/29)$ to $(35/29,0)=\mathcal L^2(15/29,-20/29)$, i.e.~$\mathcal L(\ell_1)$.\\

We need to solve $\mathcal L^4(x,y)=(x,y)$ for $(x,y) \in \mathbb{R}^2$. Note that this calculation is not trivial since $\mathcal L^4$ has $2^4=16$ affine domains to check. We summarize these computations below. Let:
\begin{eqnarray}
\mC &=& 1-1.4|x| + 0.4y, \nonumber\\
\mB &=& 1-1.4|\mC| + 0.4x, \nonumber\\
\mA &=& 1-1.4|\mB| + 0.4(\mC). \nonumber
\end{eqnarray}

Note that we need to solve,
\[\mathcal L^4 \left( \begin{array}{ccc} x \\ y \end{array} \right) = \left( \begin{array}{ccc} 1-1.4|\mA|+0.4(\mB) \\ \mA \end{array} \right)=\left( \begin{array}{ccc} x \\ y \end{array} \right).
\]
\textbf{Domain 1 and 2}:$\hspace{5mm} \mB \geq 0$, $\mC \geq 0$, $x \geq0$:
First, let us use the equality of the $y$-coordinate of $\mathcal L^4$ to $y$:
\begin{eqnarray}{\label{dom12}}
\mA = y &\implies& 1-1.4(1-1.4(1-1.4x+0.4y)+0.4x)+0.4(1-1.4x+0.4y) = y  \nonumber\\
&\implies& 0.056y = 1.96-3.864x \implies y=-69x+35.
\end{eqnarray}
Now, let us use the equality of the $x$-coordinate of $\mathcal L^4$ to $x$: \\
\emph{Domain 1}: Assuming also $\mA = y \geq0$: \\
$1-1.4y + 0.4(1-1.4(1-1.4x+0.4y)+0.4x) = x \implies -1.624y+0.84=0.056x$. Now, solving this equation together with Eqn.~\ref{dom12}, one gets $x=y=0.5$. This is the right fixed point of $\mathcal L_{1.4,0.4}$. \\
\emph{Domain 2}: Assuming $\mA = y <0$: \\
$1+1.4y + 0.4(1-1.4(1-1.4x+0.4y)+0.4x) = x \implies 1.176y+0.84=0.056x$. Now, solving this equation together with Eqn.~\ref{dom12}, one gets $x=15/29$, $y=-20/29$. This is the left end point of $\mathcal L(\ell_1)$. \\

\textbf{Domain 3 and 4}:$\hspace{5mm} \mB \geq 0$, $\mC \geq 0$, $x <0$:
From the equality of the $y$-coordinate of $\mathcal L^4$ to $y$:
\begin{eqnarray}{\label{dom34}}
\mA = y &\implies& 1-1.4(1-1.4(1+1.4x+0.4y)+0.4x)+0.4(1+1.4x+0.4y) = y  \nonumber\\
&\implies& 0.056y = 1.96+2.744x \implies y=49x+35.
\end{eqnarray}
Now, let us use the equality of the $x$-coordinate of $\mathcal L^4$ to $x$: \\
\emph{Domain 3}: Assuming also $\mA = y \geq0$: \\
$1-1.4y + 0.4(1-1.4(1+1.4x+0.4y)+0.4x) = x \implies -1.624y+0.84=1.624x$. Now, solving this equation together with Eqn.~\ref{dom34}, one gets $x=-20/29$, $y=35/29$. This is the left end point of $\ell_1$. \\
\emph{Domain 4}: Assuming $\mA = y <0$: \\
$1+1.4y + 0.4(1-1.4(1+1.4x+0.4y)+0.4x) = x \implies 1.176y+0.84=1.624x$. Now, solving this equation together with Eqn.~\ref{dom34}, one gets $x=-3/4$, $y=-7/4$. But note that at this point $\mC <0$, so this point is not in Domain 4 and there are no fixed points. \\

\textbf{Domain 5 and 6}:$\hspace{5mm} \mB \geq 0$, $\mC < 0$, $x \geq0$:
From the equality of the $y$-coordinate of $\mathcal L^4$ to $y$:
\begin{eqnarray}{\label{dom56}}
\mA = y &\implies& 1-1.4(1+1.4(1-1.4x+0.4y)+0.4x)+0.4(1-1.4x+0.4y) = y  \nonumber\\
&\implies& 1.624y = -1.96+1.624x \implies y=x-35/29.
\end{eqnarray}
Now, let us use the equality of the $x$-coordinate of $\mathcal L^4$ to $x$: \\
\emph{Domain 5}: Assuming also $\mA = y \geq0$: \\
$1-1.4y + 0.4(1+1.4(1-1.4x+0.4y)+0.4x) = x \implies -1.176y+1.96=1.624x$. Now, solving this equation together with Eqn.~\ref{dom56}, one gets $x=490/261 \approx 1.8773$, $y=175/261\approx 0.6704$. But note that at this point $\mB <0$, so this point is not in Domain 5 and there are no fixed points. \\
\emph{Domain 6}: Assuming $\mA = y <0$: \\
$1+1.4y + 0.4(1+1.4(1-1.4x+0.4y)+0.4x) = x \implies 1.624y+1.96=1.624x$. Now, solving this equation together with Eqn.~\ref{dom56}, one gets $x=x$. So, the part of the line segment $y=x-35/29$ that stays in Domain 6 is a line segment of fixed points of $\mathcal L^4$. Note that this line segment is $\mathcal L(\ell_1)$. \\

\textbf{Domain 7 and 8}:$\hspace{5mm} \mB \geq 0$, $\mC < 0$, $x <0$:
From the equality of the $y$-coordinate of $\mathcal L^4$ to $y$:
\begin{eqnarray}{\label{dom78}}
\mA = y &\implies& 1-1.4(1+1.4(1+1.4x+0.4y)+0.4x)+0.4(1+1.4x+0.4y) = y  \nonumber\\
&\implies& 0.84y = -1.96-2.744x \implies y=(49/15)x-7/3.
\end{eqnarray}
Now, let us use the equality of the $x$-coordinate of $\mathcal L^4$ to $x$: \\
\emph{Domain 7}: Assuming also $\mA = y \geq0$: \\
$1-1.4y + 0.4(1+1.4(1+1.4x+0.4y)+0.4x) = x \implies -1.176y+1.96=0.056x$. Now, solving this equation together with Eqn.~\ref{dom78}, one gets $x=35/29 \approx 1.2068$, $y=140/87\approx 1.6091$. But note that at this point $x\geq0$, so this point is not in Domain 7 and there are no fixed points. \\
\emph{Domain 8}: Assuming $\mA = y <0$: \\
$1+1.4y + 0.4(1+1.4(1+1.4x+0.4y)+0.4x) = x \implies 1.624y+1.96=0.056x$. Now, solving this equation together with Eqn.~\ref{dom78}, one gets $x \approx 0.3485$, $y\approx -1.1948$. But note that at this point $x\geq0$, so this point is not in Domain 8 and there are no fixed points. \\

\textbf{Domain 9 and 10}:$\hspace{5mm} \mB < 0$, $\mC \geq 0$, $x \geq0$:
From the equality of the $y$-coordinate of $\mathcal L^4$ to $y$:
\begin{eqnarray}{\label{dom910}}
\mA = y &\implies& 1+1.4(1-1.4(1-1.4x+0.4y)+0.4x)+0.4(1-1.4x+0.4y) = y  \nonumber\\
&\implies& 1.624y = 0.84+2.744x \implies y=(49/29)x+15/29.
\end{eqnarray}
Now, let us use the equality of the $x$-coordinate of $\mathcal L^4$ to $x$: \\
\emph{Domain 9}: Assuming also $\mA = y \geq0$: \\
$1-1.4y + 0.4(1-1.4(1-1.4x+0.4y)+0.4x) = x \implies -1.624y+0.84=0.056x$. Now, solving this equation together with Eqn.~\ref{dom910}, one gets $x=0$, $y=15/29$. This is the right end point of $\ell_1$. \\
\emph{Domain 10}: Assuming $\mA = y <0$: \\
$1+1.4y + 0.4(1-1.4(1-1.4x+0.4y)+0.4x) = x \implies 1.176y+0.84=0.056x$. Now, solving this equation together with Eqn.~\ref{dom910}, one gets $x=y=-3/4$. But note that at this point $x<0$, so this point is not in Domain 10 and there are no fixed points. \\

\textbf{Domain 11 and 12}:$\hspace{5mm} \mB < 0$, $\mC \geq 0$, $x <0$:
From the equality of the $y$-coordinate of $\mathcal L^4$ to $y$:
\begin{eqnarray}{\label{dom1112}}
\mA = y &\implies& 1+1.4(1-1.4(1+1.4x+0.4y)+0.4x)+0.4(1+1.4x+0.4y) = y  \nonumber\\
&\implies& 1.624y = 0.84-1.624x \implies y=-x+15/29.
\end{eqnarray}
Now, let us use the equality of the $x$-coordinate of $\mathcal L^4$ to $x$: \\
\emph{Domain 11}: Assuming also $\mA = y \geq0$: \\
$1-1.4y + 0.4(1-1.4(1+1.4x+0.4y)+0.4x) = x \implies -1.624y+0.84=1.624x$. Now, solving this equation together with Eqn.~\ref{dom1112}, one gets $x=x$. So, the part of the line segment $y=-x-15/29$ that stays in Domain 11 is a line segment of fixed points of $\mathcal L^4$. Note that this line segment is $\mathcal L(\ell_1)$. \\
\emph{Domain 12}: Assuming $\mA = y <0$: \\
$1+1.4y + 0.4(1-1.4(1+1.4x+0.4y)+0.4x) = x \implies 1.176y+0.84=1.624x$. Now, solving this equation together with Eqn.~\ref{dom1112}, one gets $x\approx3.8813$, $y\approx-3.3641$. But note that at this point $x\geq0$, so this point is not in Domain 12 and there are no fixed points. \\

\textbf{Domain 13 and 14}:$\hspace{5mm} \mB < 0$, $\mC < 0$, $x \geq0$:
From the equality of the $y$-coordinate of $\mathcal L^4$ to $y$:
\begin{eqnarray}{\label{dom1314}}
\mA = y &\implies& 1+1.4(1+1.4(1-1.4x+0.4y)+0.4x)+0.4(1-1.4x+0.4y) = y  \nonumber\\
&\implies& 0.056y = 4.76-2.744x \implies y=-49x+85.
\end{eqnarray}
Now, let us use the equality of the $x$-coordinate of $\mathcal L^4$ to $x$: \\
\emph{Domain 13}: Assuming also $\mA = y \geq0$: \\
$1-1.4y + 0.4(1+1.4(1-1.4x+0.4y)+0.4x) = x \implies -1.176y+1.96=1.624x$. Now, solving this equation together with Eqn.~\ref{dom1314}, one gets $x=1.75$, $y=-0.75$. But note that at this point $y=\mA <0$ so this point is not in Domain 13 and there are no fixed points. \\
\emph{Domain 14}: Assuming $\mA = y <0$: \\
$1+1.4y + 0.4(1+1.4(1-1.4x+0.4y)+0.4x) = x \implies 1.624y+1.96=1.624x$. Now, solving this equation together with Eqn.~\ref{dom1314}, one gets $x\approx1.7241$, $y\approx0.5172$. But note that at this point $y=\mA \geq0$, so this point is not in Domain 14 and there are no fixed points. \\

\textbf{Domain 15 and 16}:$\hspace{5mm} \mB < 0$, $\mC < 0$, $x < 0$:
From the equality of the $y$-coordinate of $\mathcal L^4$ to $y$:
\begin{eqnarray}{\label{dom1516}}
\mA = y &\implies& 1+1.4(1+1.4(1+1.4x+0.4y)+0.4x)+0.4(1+1.4x+0.4y) = y  \nonumber\\
&\implies& 0.056y = 4.76+3.864x \implies y=69x+85.
\end{eqnarray}
Now, let us use the equality of the $x$-coordinate of $\mathcal L^4$ to $x$: \\
\emph{Domain 15}: Assuming also $\mA = y \geq0$: \\
$1-1.4y + 0.4(1+1.4(1+1.4x+0.4y)+0.4x) = x \implies -1.176y+1.96=0.056x$. Now, solving this equation together with Eqn.~\ref{dom1516}, one gets $x=-35/29$, $y=50/29$. But note that at this point $\mC \geq0$, so this point is not in Domain 15 and there are no fixed points. \\
\emph{Domain 16}: Assuming $\mA = y <0$: \\
$1+1.4y + 0.4(1+1.4(1+1.4x+0.4y)+0.4x) = x \implies 1.624y+1.96=0.056x$. Now, solving this equation together with Eqn.~\ref{dom1516}, one gets $x=y=-5/4$. This is the left fixed point of $\mathcal L_{1.4,0.4}$. \\

\paragraph{Acknowledgments}
I would like to thank S.~E.~Newhouse for his helpful discussions and suggestions. I also would like to thank Duncan Sands for providing corrections to some historical comments and the anonymous referee whose comments improved the exposition of the paper.

\bibliography{references}
\bibliographystyle{plain}

\ttfamily
Department of Mathematics, Michigan State University, East Lansing, MI 48824 USA \\
Current address: Max Planck Institute for Cognitive and Brain Sciences, Leipzig, Germany \\ E-mail: yildiz@cbs.mpg.de
\end{document}

%% file: macros.tex
\newcommand{\ind}{\mbox{\mbox{ } \mbox{ } \mbox{ }}}

\newcommand{\be}{\begin{enumerate}}
\newcommand{\ee}{\end{enumerate}}
\newcommand{\bit}{\begin{itemize}}
\newcommand{\eit}{\end{itemize}}
\newcommand{\bc}{\begin{center}}
\newcommand{\ec}{\end{center}}
\newcommand{\beq}{\begin{equation}}
\newcommand{\eeq}{\end{equation}}
\newcommand{\bqn}{\begin{eqnarray}}
\newcommand{\eqn}{\end{eqnarray}}
\newcommand{\bqns}{\begin{eqnarray*}}
\newcommand{\eqns}{\end{eqnarray*}}
\newcommand{\N}{{\bf N}}

\newcommand{\R}{{\bf R}}

\newcommand{\rarrow}{\rightarrow}

\newcommand{\mA}{{\mathcal A}}
\newcommand{\mB}{{\mathcal B}}
\newcommand{\mC}{{\mathcal C}}

\newcommand{\se}{{\mathcal E}}

\newcommand{\bdoc}{\begin{document}}
\newcommand{\edoc}{\end{document}}
\newcommand{\docart}{\documentstyle[12pt]{article}}
\newcommand{\vs}{\vspace{.16in}}
\newcommand{\half}{\frac{1}{2}}
\newcommand{\third}{\frac{1}{3}}
\newcommand{\fourth}{\frac{1}{4}}
\newcommand{\diffeo}{\mbox{ diffeomorphism }}
\newcommand{\dnorm}[1]{\mbox{$\mid \mid #1 \mid \mid$}}
\newcommand{\zero}{ \mbox{$ {\bf 0 } $}}
\newcommand{\PL}{\mbox{$Per_L(\R)$}}

\newcommand{\0}{{\bf 0}}

\newcommand{\twoover}[2]{{\footnotesize \begin{array}[t]{c}\mbox{ {\normalsize #1}}  \\
                               #2
        \end{array} } \hspace{1cm} }
\newcommand{\threeover}[3]{ {\footnotesize \begin{array}[t]{c}\mbox{ {\normalsize #1}}  \\
                               #2  \\
                                 #3
        \end{array} } \hspace{1cm} }
\newcommand{\fourover}[4]{ {\footnotesize \begin{array}[t]{c}\mbox{ {\normalsize #1}}  \\
                               #2  \\
                                 #3 \\
				#4
        \end{array} } \hspace{1cm}  }

\newenvironment{vover3}[3]{
       {\footnotesize \begin{array}[t]{c}
                               {\normalsize #1} \vspace{.2cm} \\
                                 #2 \\
                                 #3
                       \end{array} }
                            } { \hspace{1cm} }
\newenvironment{vover4}[4] {
       {\footnotesize \begin{array}[t]{c}
                               {\normalsize #1} \vspace{.2cm} \\
                                 #2 \\
                                 #3 \\
                                 #4
                       \end{array} }
                            } {  \hspace{1cm} }

\newenvironment{seteqover}[3]{#1 \ \ = \hspace{.6cm}
                                #2 #3} {\ind}

\newcommand{\Norm}[1]{\mbox{$\parallel #1 \parallel$}}
\newcommand{\imp}{\mbox{$\Longrightarrow$}}
\newcommand{\pr}{\mbox{$\prime$}}

\newcommand{\twocvec}[2]{\mbox{$\left( \begin{array}{c}
				 #1 \\
				 #2
				\end{array} \right)$} }

\newcommand{\twomat}[4]{\mbox{$ \left( \begin{array}{cc}
                              #1 & #2 \\
                           #3 & #4
                       \end{array} \right) $ } }

\newcommand{\Itwomat}[4]{\mbox{$ \left( \begin{array}{cc}
                              #1 & #2 \\
                           #3 & #4
                       \end{array} \right)^{-1} $ } }

\newcommand{\threemat}[9]{ \[ \left( \begin{array}{ccc}
				 #1 & #2 & #3 \\
				 #4 & #5 & #6 \\
				 #7 & #8 & #9
				\end{array} \right) \]  }

\newcommand{\newthreemat}[9]{ \left( \begin{array}{ccc}
				 #1 & #2 & #3 \\
				 #4 & #5 & #6 \\
				 #7 & #8 & #9
				\end{array} \right) }

\newcommand{\Ithreemat}[9]{\mbox{ $\left( \begin{array}{ccc}
				 #1 & #2 & #3 \\
				 #4 & #5 & #6 \\
				 #7 & #8 & #9
				\end{array} \right)^{-1} $ } }

\newcommand{\ba}{\mbox{${\bf a}$}}
\newcommand{\bfa}{\mbox{${\bf a}$}}
\newcommand{\bb}{\mbox{${\bf b}$}}
\newcommand{\bfb}{\mbox{${\bf b}$}}
\newcommand{\bfc}{\mbox{${\bf c}$}}
\newcommand{\bi}{\mbox{${\bf i}$}}
\newcommand{\bj}{\mbox{${\bf j}$}}
\newcommand{\bfi}{\mbox{${\bf i}$}}
\newcommand{\bfj}{\mbox{${\bf j}$}}
\newcommand{\bfk}{\mbox{${\bf k}$}}
\newcommand{\bfx}{\mbox{${\bf x}$}}
\newcommand{\bfy}{\mbox{${\bf y}$}}
\newcommand{\bfz}{\mbox{${\bf z}$}}
\newcommand{\bfu}{\mbox{${\bf u}$}}
\newcommand{\bfv}{\mbox{${\bf v}$}}
\newcommand{\bfw}{\mbox{${\bf w}$}}
\newcommand{\bfg}{\mbox{${\bf g}$}}
\newcommand{\onef}{\mbox{${\bf 1}$}}
\renewcommand{\bfc}{\mbox{${\bf c}$}}
\newcommand{\bfe}{\mbox{${\bf e}$}}
\newcommand{\bff}{\mbox{${\bf f}$}}
\newcommand{\bS}{\mbox{${\bf S}$}}
\newcommand{\bV}{\mbox{${\bf V}$}}
\newcommand{\bw}{\bar{w}}
\newcommand{\fh}{\hat{f}}
\newcommand{\wh}{\hat{w}}
\newcommand{\Wh}{\hat{W}}
\newcommand{\Ih}{\hat{I}}
\newcommand{\Lh}{\hat{L}}
\newcommand{\grgh}{\hat{\gamma}}
\newcommand{\br}{\mbox{${\bf r}$}}
\newcommand{\bs}{\mbox{${\bf s}$}}

\newcommand{\oz}{\overline{z}}
\newcommand{\ogz}{\overline{gz}}
\newcommand{\oxi}{\overline{\xi}}

\newcommand{\ninf}{n \rarrow \infty}
\newcommand{\linf}{\lim_{n \rarrow \infty}}
\newcommand{\fon}{ \frac{1}{n} }
\newcommand{\skn}{ \sum_{k=0}^{n-1} }
\newcommand{\siin}{ \sum_{i=0}^{n-1} }

\newcommand{\alp}{\mbox{$\alpha$}}
\newcommand{\gra}{\mbox{$\alpha$}}

\newcommand{\bet}{\mbox{$\beta$}}
\newcommand{\grb}{\mbox{$\beta$}}

\newcommand{\gam}{\mbox{$\gamma$}}
\newcommand{\Gam}{\mbox{$\Gamma$}}
\newcommand{\Del}{\mbox{$\Delta$}}
\newcommand{\sig}{\mbox{$\sigma$}}
\newcommand{\Sig}{\mbox{$\Sigma$}}

\newcommand{\grg}{\mbox{$\gamma$}}
\newcommand{\Grg}{\mbox{$\Gamma$}}
\newcommand{\grd}{\mbox{$\delta$}}
\newcommand{\Grd}{\mbox{$\Delta$}}
\newcommand{\grs}{\mbox{$\sigma$}}
\newcommand{\Grs}{\mbox{$\Sigma$}}
\newcommand{\grf}{\mbox{$\phi$}}
\newcommand{\Grf}{\mbox{$\Phi$}}
\newcommand{\grth}{\mbox{$\theta$}}
\newcommand{\Grth}{\mbox{$\Theta$}}
\newcommand{\gro}{\mbox{$\omega$}}
\newcommand{\Gro}{\mbox{$\Omega$}}
\newcommand{\gre}{\mbox{$\epsilon$}}
\newcommand{\bgre}{\mbox{$\bar{\epsilon}$}}
\newcommand{\grve}{\mbox{$\varepsilon$}}
\newcommand{\grr}{\mbox{$\rho$}}
\newcommand{\Grr}{\mbox{$\Rho$}}
\newcommand{\grt}{\mbox{$\tau$}}
\newcommand{\gri}{\mbox{$\iota$}}
\newcommand{\grp}{\mbox{$\pi$}}
\newcommand{\grP}{\mbox{$\Pi$}}
\newcommand{\grl}{\mbox{$\lambda$}}
\newcommand{\grlh}{\mbox{$\hat{\lambda}$}}
\newcommand{\grL}{\mbox{$\Lambda$}}
\newcommand{\Grl}{\mbox{$\Lambda$}}
\newcommand{\grz}{\mbox{$\zeta$}}

\newcommand{\thet}{\mbox{$\theta$}}
\newcommand{\Thet}{\mbox{$\Theta$}}
\newcommand{\ome}{\mbox{$\omega$}}
\newcommand{\Ome}{\mbox{$\Omega$}}
\newcommand{\eps}{\mbox{$\epsilon$}}

\newcommand{\iot}{\mbox{$\iota$}}

\newcommand{\lam}{\mbox{$\lambda$}}
\newcommand{\Lam}{\mbox{$\Lambda$}}
\newcommand{\zet}{\mbox{$\zeta$}}

\newcommand{\kap}{\mbox{$\kappa$}}

\newcommand{\grx}{\mbox{$\xi$}}
\newcommand{\Grx}{\mbox{$\Xi$}}
\newcommand{\grc}{\mbox{$\chi$}}
\newcommand{\Grc}{\mbox{$\Chi$}}
\newcommand{\grn}{\mbox{$\nu$}}
\newcommand{\grm}{\mbox{$\mu$}}
\newcommand{\grk}{\mbox{$\kappa$}}

\newtheorem{theorem}{Theorem}[section]
\newtheorem{Cor}[theorem]{Corollary}
\newtheorem{Corollary}[theorem]{Corollary}
\newtheorem{Prop}[theorem]{Proposition}
\newtheorem{Proposition}[theorem]{Proposition}
\newtheorem{Definition}[theorem]{Definition}
\newtheorem{Lem}[theorem]{Lemma}
\newtheorem{Lemma}[theorem]{Lemma}
\newtheorem{Rem}[theorem]{Remark}
\newtheorem{Remark}[theorem]{Remark}
\newtheorem{Conjecture}[theorem]{Conjecture}
\newtheorem{Theorem}[theorem]{Theorem}
\newtheorem{Th}[theorem]{Theorem}
\newtheorem{question}{Question}
\newtheorem{conjecture}{Conjecture}

\newcommand{\ta}{\tilde{a}}
\newcommand{\tgra}{\tilde{\alpha}}
\newcommand{\tGrl}{\tilde{\Grl}}
\newcommand{\tb}{\tilde{b}}
\newcommand{\tgrb}{\tilde{\beta}}
\newcommand{\tgrg}{\tilde{\gamma}}
\newcommand{\tc}{\tilde{c}}
\newcommand{\td}{\tilde{d}}
\newcommand{\te}{\tilde{e}}
\newcommand{\tf}{\tilde{f}}
\newcommand{\teta}{\tilde{\eta}}
\newcommand{\tg}{\tilde{g}}
\newcommand{\tih}{\tilde{h}}
\renewcommand{\th}{\tilde{h}}
\newcommand{\ti}{\tilde{i}}
\newcommand{\tj}{\tilde{j}}
\newcommand{\tk}{\tilde{k}}
\newcommand{\tl}{\tilde{l}}
\newcommand{\tm}{\tilde{m}}
\newcommand{\tn}{\tilde{n}}
\newcommand{\tio}{\tilde{o}}
\newcommand{\tp}{\tilde{p}}
\newcommand{\tq}{\tilde{q}}
\newcommand{\tr}{\tilde{r}}
\newcommand{\ts}{\tilde{s}}
\newcommand{\tit}{\tilde{t}}
\newcommand{\tu}{\tilde{u}}
\newcommand{\tv}{\tilde{v}}
\newcommand{\tw}{\tilde{w}}
\newcommand{\tx}{\tilde{x}}
\newcommand{\ty}{\tilde{y}}
\newcommand{\tz}{\tilde{z}}
\newcommand{\tA}{\tilde{A}}
\newcommand{\tB}{\tilde{B}}
\newcommand{\tC}{\tilde{C}}
\newcommand{\tD}{\tilde{D}}
\newcommand{\tE}{\tilde{E}}
\newcommand{\tF}{\tilde{F}}
\newcommand{\tG}{\tilde{G}}
\newcommand{\tH}{\tilde{H}}
\newcommand{\tI}{\tilde{I}}
\newcommand{\tJ}{\tilde{J}}
\newcommand{\tK}{\tilde{K}}
\newcommand{\tL}{\tilde{L}}
\newcommand{\tM}{\tilde{M}}
\newcommand{\tN}{\tilde{N}}
\newcommand{\tO}{\tilde{O}}
\newcommand{\tP}{\tilde{P}}
\newcommand{\tQ}{\tilde{Q}}
\newcommand{\tR}{\tilde{R}}
\newcommand{\tS}{\tilde{S}}
\newcommand{\tT}{\tilde{T}}
\newcommand{\tU}{\tilde{U}}
\newcommand{\tV}{\tilde{V}}
\newcommand{\tW}{\tilde{W}}
\newcommand{\tX}{\tilde{X}}
\newcommand{\tY}{\tilde{Y}}
\newcommand{\tZ}{\tilde{Z}}

\newcommand{\tTh}{\tilde{\Theta}}

\newcommand{\fl}{f_{\grl}}
\newcommand{\el}{E_{\grl}}
\newcommand{\sel}{\se_{\grl}}
\newcommand{\szl}{\sz_{\grl}}
\newcommand{\zl}{Z_{\grl}}
\newcommand{\famfl}{\{ \fl \} }
\newcommand{\sT}{{\mathcal T}}
\newcommand{\hsT}{{\hat{\mathcal T}}}
\newcommand{\hsH}{{\hat{\mathcal H}}}
\newcommand{\hsCZ}{{\hat{\mathcal CZ}}}
\newcommand{\sCG}{{\mathcal CG}}
\newcommand{\sFG}{{\mathcal FG}}
\newcommand{\sTemp}{{\mathcal TM}}
\newcommand{\ab}{{\bf a}}
\newcommand{\xb}{{\bf x}}
\newcommand{\yb}{{\bf y}}
\newcommand{\Def}{{\bf Definition}}
\newcommand{\Pf}{{\bf Proof}}
\newcommand{\uG}{\bar{{\mathcal G}}}
\newcommand{\sG}{{\mathcal G}}
\newcommand{\sF}{{\mathcal F}}
\newcommand{\stG}{\tilde{{\mathcal G}}}
\newcommand{\tmu}{\tilde{\mu}}
\newcommand{\tpi}{\tilde{\pi}}
\newcommand{\tgrl}{\tilde{\grl}}
\newcommand{\sD}{{\mathcal D}}
\newcommand{\sDh}{\hat{\sD}}
\newcommand{\sP}{{\mathcal P}}
\newcommand{\tsP}{{\tilde{\mathcal P}}}
\newcommand{\sPh}{\hat{\sP}}
\newcommand{\sM}{{\mathcal M}}
\newcommand{\bsM}{\bar{{\mathcal M}}}
\newcommand{\sN}{{\mathcal N}}
\newcommand{\sS}{{\mathcal S}}
\newcommand{\tsS}{\tilde{{\mathcal S}}}
\newcommand{\tgrt}{\tilde{{\tau}}}
\newcommand{\sI}{{\mathcal I}}
\newcommand{\sW}{{\mathcal W}}
\newcommand{\sCZ}{{\mathcal C}{\mathcal Z}}
\newcommand{\tsT}{\tilde{\sT}}
\newcommand{\tsC}{\tilde{\sC}}
\newcommand{\tsD}{\tilde{\sD}}
\newcommand{\tgre}{\tilde{\gre}}
\newcommand{\sH}{{\mathcal H}}
\newcommand{\tsH}{\tilde{\sH}}
\newcommand{\tsG}{\tilde{\sG}}
\newcommand{\tsW}{\tilde{\sW}}
\newcommand{\tDel}{\tilde{\Grd}}
\newcommand{\sA}{{\mathcal A}}
\newcommand{\sB}{{\mathcal B}}
\newcommand{\sCF}{{\mathcal CF}}
\newcommand{\bsH}{\bar{\sH}}
\newcommand{\ntwob}{\bar{\nu_2}}
\newcommand{\nb}{\bar{\nu}}
\newcommand{\bE}{\bar{E}}
\newcommand{\bF}{\bar{F}}
\newcommand{\bG}{\bar{G}}
\newcommand{\bI}{\bar{I}}
\newcommand{\bZ}{\bar{Z}}
\newcommand{\bh}{\bar{h}}
\newcommand{\bm}{\bar{m}}
\newcommand{\ml}{m_{\ell}}
\newcommand{\bT}{\bar{T}}
\newcommand{\bK}{\bar{K}}
\newcommand{\bx}{\bar{x}}
\newcommand{\bB}{\bar{B}}
\newcommand{\bxi}{\bar{\xi}}
\newcommand{\bu}{\bar{u}}
\newcommand{\bv}{\bar{v}}
\newcommand{\bg}{\bar{g}}
\newcommand{\bN}{\bar{\N}}
\newcommand{\bX}{\bar{X}}
\newcommand{\bY}{\bar{Y}}
\newcommand{\bGrs}{\bar{\Grs}}
\newcommand{\bgrs}{\bar{\grs}}
\newcommand{\bmu}{\bar{\mu}}
\newcommand{\bnu}{\bar{\nu}}
\newcommand{\bTh}{\bar{\Theta}}
\newcommand{\db}{\bar{d}}
\newcommand{\dx}{\dot{x}}
\newcommand{\dy}{\dot{y}}
\newcommand{\dz}{\dot{z}}
\newcommand{\ddx}{\partial_x}
\newcommand{\ddy}{\partial_y}
\newcommand{\ddz}{\partial_z}
\newcommand{\dxn}{\dot{x}_n}
\newcommand{\n}{\norm}
\newcommand{\No}{\dnorm}
\newcommand{\no}[1]{\mbox{$\norm{#1}_1$}}
\newcommand{\grfh}{\hat{\grf}}

\newcommand{\ol}{\overline}
\newcommand{\ob}{\bar}

\newcommand{\eqdef}{\mbox{ $\stackrel{{\rm def}}{=}$ }}
\newcommand{\defeq}{\mbox{ $\stackrel{{\rm def}}{=}$ }}

\newcommand{\ueps}{\underline{\grve}}
\newcommand{\uepss}{\underline{\grve}^s}
\newcommand{\uepsu}{\underline{\grve}^u}
\newcommand{\headeps}{(\varepsilon_{0},\varepsilon_{1} \ldots )}
\newcommand{\taileps}{(\ldots, \varepsilon_{-2},\varepsilon_{-1})}
\newcommand{\fulleps}{(\ldots \varepsilon_{-2},\varepsilon_{-1}\cdot\varepsilon_{0},\varepsilon_{1}\ldots) }

\newcommand{\lex}[2]{{\vphantom{#2}}_{#1}{#2}} 

\newcommand{\Lpfk}{\mbox{$\cap \hspace{-.44cm}\mid$}}

%% file: Discontinuity_yildiz_etds.bbl
\begin{thebibliography}{10}

\bibitem{Buzzi}
J.~Buzzi.
\newblock Maximal entropy measures for piecewise affine surface homeomorphisms.
\newblock {\em Ergod. th. dynam. systems}, 29:1723--1763, 2009.

\bibitem{Zyg3}
Z.~Galias.
\newblock Obtaining rigorous bounds for topological entropy for discrete time
  dynamical systems.
\newblock {\em Proc. Internat. Symposium on Nonlinear Theory and its
  Applications}, pages 619 -- 622, 2002.

\bibitem{Zyg2}
Z.~Galias and P.~Zygliczy\'{n}ski.
\newblock Abundance of homoclinic and heteroclinic orbits and rigorous bounds
  for the topological entropy for the {H}\'enon map.
\newblock {\em Nonlinearity}, 14:909 -- 932, 2001.

\bibitem{Henon}
M.~H\'{e}non.
\newblock A two-dimensional mapping with a strange attractor.
\newblock {\em Communications in Mathematical Physics}, 50:69--77, 1976.

\bibitem{Ikeda}
K.~Ikeda, H.~Daido, and Akimoto O.
\newblock Optical turbulence: {C}haotic behavior of transmitted light from a
  ring cavity.
\newblock {\em Phys. Rev. Lett.}, 45:709--712, 1980.

\bibitem{Ish3}
Y.~Ishii and D.~Sands.
\newblock Monotonicity of the {L}ozi {F}amily {N}ear the {T}ent-{M}aps.
\newblock {\em Commun. Math. Phys.}, 198:397--406, 1998.

\bibitem{Ish4}
Y.~Ishii and D.~Sands.
\newblock Lap number entropy formula for piecewise affine and projective maps
  in several dimensions.
\newblock {\em Nonlinearity}, 20:2755--2772(18), 2007.

\bibitem{Katok1}
A.~Katok.
\newblock Lyapunov exponents, entropy and periodic orbits for diffeomorphisms.
\newblock {\em Inst. Hautes Etudes Sci. Publ. Math.}, 51(1):137--173, 1980.

\bibitem{Katok2}
A.~Katok.
\newblock Nonuniform hyperbolicity and structure of smooth dynamical systems.
\newblock {\em Proc. of Intl. Congress of Math.}, 2:1245--1254, 1983.

\bibitem{Mis4}
M.~Misiurewicz.
\newblock On non-continuity of topological entropy.
\newblock {\em Bull. Acad. Polon. Sci., Ser. Sci. Math. Astro. Phys.},
  19(4):319--320, 1971.

\bibitem{Mis2}
M.~Misiurewicz.
\newblock Diffeomorphisms without any measure with maximal entropy.
\newblock {\em Bull. Acad. Polon. Sci., Ser. Sci. Math. Astro. Phys.},
  21(10):903--910, 1973.

\bibitem{Mis5}
M.~Misiurewicz.
\newblock Jumps of entropy in one dimension.
\newblock {\em Fund. Math.}, 132(3):215--226, 1989.

\bibitem{Mis1}
M.~Misiurewicz and W.~Szlenk.
\newblock Entropy of piecewise monotone mappings.
\newblock {\em Studia Mathematica}, 67(1):45--63, 1980.

\bibitem{Newhouse2}
S.~Newhouse.
\newblock Continuity properties of entropy.
\newblock {\em Ann. of Math.}, 129:215--235, 1989.

\bibitem{Newhouse}
S.~Newhouse, M.~Berz, J.~Grote, and K.~Makino.
\newblock On the estimation of topological entropy on surfaces.
\newblock {\em Contemporary Mathematics}, 469:243--270, 2008.

\bibitem{Rees}
M.~Rees.
\newblock A minimal positive entropy homeomorphism of the 2-torus.
\newblock {\em J. London Math. Soc. (2)}, 23:537--550, 1981.

\bibitem{Burak}
I.B. Yildiz.
\newblock Monotonicity of the {L}ozi family and the zero entropy locus.
\newblock {\em Nonlinearity}, 24:1613--1628, 2011.

\bibitem{Yomdin}
Y.~Yomdin.
\newblock Volume growth and entropy.
\newblock {\em Israel J. Math.}, 57:285--300, 1987.

\bibitem{Zyg}
P.~Zygliczy\'{n}ski.
\newblock Computer assisted proof of chaos in the {R}ossler equations and the
  {H}\'{e}non map.
\newblock {\em Nonlinearity}, 10(1):243 -- 252, 1997.

\end{thebibliography}
